\documentclass[]{amsart}
\usepackage{amsmath}
\usepackage{amssymb}
\usepackage{latexsym,boxedminipage}
\newcommand{\la}{\lambda}

\newtheorem{theorem}{Theorem}
\newtheorem{corollary}{Corollary}                                 
                                 
\newtheorem{proposition}{Proposition}                                 
\title{Lecture Hall Theorems, $q$-series and Truncated Objects}
\date{}
\author{Sylvie Corteel}
\thanks{Supported by ATIP Jeune chercheur CNRS}
\address{CNRS PRiSM, UVSQ, Versailles, France}
\email{sylvie.corteel@prism.uvsq.fr}
\author{Carla D. Savage}
\thanks{Research supported by NSA grant MDA 904-01-0-0083 and NSF INT-0230800}
\address{Dept. of Computer Science,
N. C. State University, Raleigh, USA}
\email{savage@csc.ncsu.edu}

\begin{document}
\maketitle
\begin{abstract}
We show here that the refined theorems for both lecture hall
partitions and anti-lecture hall compositions can be obtained as
straightforward consequences of 
two $q$-Chu Vandermonde identities, once an appropriate recurrence is
derived.
%two $q$-Chu Vandermonde identities~:
%$$
%\sum_{m=0}^n \frac{(a;q)_m(q^{-n};q)_m}
%{(c;q)_m(q;q)_m}q^m=\frac{a^n(c/a;q)_n}{(c;q)_n}.
%$$            
%and 
%$$\sum_{m=0}^n \left[ \begin{matrix} n \\ m
%\end{matrix} \right]_q    \frac{(a;q)_m}{(c;q)_m}(-c/a)^mq^{m\choose 2}
%=\frac{(c/a;q)_n}{(c;q)_n}.
%$$   
We use this approach to get new lecture hall-type theorems
for truncated objects. The {\em truncated lecture hall partitions}
are sequences $(\lambda_1,\ldots ,\lambda_k)$
such that
\[
\frac{\lambda_1}{n}\ge \frac{\lambda_2}{n-1}\ge \ldots \ge 
\frac{\lambda_k}{n-k+1}\ge 0
\]
and we show that their generating function is~:
\[
\sum_{m=0}^k \left[ \begin{matrix} n \\ m
\end{matrix} \right]_q
q^{m+1\choose 2}
  \frac{(-q^{n-m+1};q)_m}{(q^{2n-m+1};q)_m}.
\] From this,
we are able to give a combinatorial characterization of
truncated lecture hall partitions and new finitizations of
refinements of Euler's theorem. 
The {\em truncated anti-lecture hall compositions}
are sequences $(\lambda_1,\ldots ,\lambda_k)$
such that
\[
\frac{\lambda_1}{n-k+1}\ge \frac{\lambda_2}{n-k+2}\ge \ldots \ge 
\frac{\lambda_k}{n}\ge 0.
\]
We show that their generating function is~:
\[
 \left[ \begin{matrix} n \\ k
\end{matrix} \right]_q
  \frac{(-q^{n-k+1};q)_k}{(q^{2(n-k+1)};q)_k},
\]
%\[ 
%\frac{(q^{2(n-k+1)};q^2)_k}{(q;q)_k(q^{2(n-k+1)};q)_k},
%\]
giving a finitization of a well-known partition identity.
We give two different multivariate refinements of these new results~:
the $q$-calculus approach gives $(u,v,q)$-refinements,
 while a completely different approach gives odd/even
$(x,y)$-refinements.

\end{abstract}

\section{Introduction}

For a sequence $\lambda = (\lambda_1, \lambda_2, \ldots, \lambda_n)$
of nonnegative integers,
define the {\em weight} of $\lambda$ to be
$|\lambda|=\lambda_1 + \cdots + \lambda_n$ and call each $\lambda_i$
a {\em part} of $\lambda$.  If  $\lambda$ 
has all parts nonnegative, we call it a {\em composition} 
 and if, in addition, $\lambda$ is a
nonincreasing sequence, we call
it a {\em partition}.

In \cite{BME1}, inspired by work of Eriksson and Eriksson on Coxeter groups,
Bousquet-M{\'e}lou and Eriksson considered 
{\em lecture hall partitions}, specifically, the set $L_n$
of partitions, $\lambda$, into  $n$ nonnegative parts satisfying

\[
\frac{\la_1}{n}\ge\frac{\la_{2}}{n-1}\ge \ldots \ge \frac{\la_{n-1}}{2} 
\geq
 \frac{\la_n}{1}\ge 0,
\]
and proved the following surprising result.

\noindent
{\bf The  Lecture Hall Theorem \cite{BME1}~:}
\begin{equation}
L_n(q) \triangleq  \sum_{\lambda \in L_n} q^{|\lambda|} = \frac{1}{(q;q^2)_n}
\label{lhp}
\end{equation}
with $(a;q)_n=\prod_{i=0}^{n-1}(1-aq^i)$.

In \cite{BME1}, Bousquet-M\'elou and Eriksson gave two proofs of the
Lecture Hall Theorem~: one
based on Coxeter groups and one combinatorial.
Subsequently, Andrews gave a  proof based on partition analysis
\cite{Om1}. A refinement and generalizations
of the identity (\ref{lhp}) were given by Bousquet-M\'elou and Eriksson
in \cite{BME2}. They showed that an elegant mapping between certain partitions
in $L_n$ and partitions in $L_{n-1}$ gives a functional equation
which easily implies the result.  
The first bijective proof of
the Lecture Hall Theorem was given by Yee
\cite{yee1} and is close to \cite{BME1}. Others followed \cite{Om5, ne}.

An involution
in \cite{BME1} gave the following nice refinement 
of (\ref{lhp}).

\noindent
{\rm \bf The Odd/Even Lecture Hall Theorem  \cite{BME1}~:}
Given $\lambda=(\lambda_1,\lambda_2,\ldots)$, define $\lambda_o=
(\lambda_1,\lambda_3,\ldots)$ and $\lambda_e=
(\lambda_2,\lambda_4,\ldots)$. Then
\begin{equation}
L_n(x,y)  \triangleq
\sum_{\lambda \in L_n} x^{|\lambda_o|} y^{|\lambda_e|} = \prod_{i=1}^n
\frac{1}{1-x^iy^{i-1}}.
\label{lhpxy}
\end{equation}

A different approach of Bousquet-M{\'e}lou and Eriksson
in \cite{BME3} led to the following.

\noindent
{\rm \bf  The Refined Lecture Hall Theorem \cite{BME3}~:}
\begin{equation}
L_n(u,v,q) \triangleq \sum_{\lambda \in L_n} q^{|\lambda|}
u^{|\lceil \lambda \rceil |}
v^{o(\lceil \lambda \rceil)}  = 
\frac{(-uvq;q)_n}{(u^2q^{n+1};q)_n},
\label{lhpuv}
\end{equation}
where for a partition $\la = (\la_1, \ldots, \la_n)$,
$o(\la)$ is the number of odd parts of $\la$ and
$\lceil \la \rceil$ is the partition
$(\lceil \la_1/n \rceil,
\lceil \la_2/(n-1) \rceil, \ldots,
\lceil \la_{n-1}/2 \rceil, 
\lceil \la_n/1 \rceil)$.

Setting $u=v=1$ in (\ref{lhpuv}) gives (\ref{lhp}).
Yee gave  a beautiful bijective proof \cite{yee2} of
this theorem.

In \cite{CS}, we considered a new twist on these results
by studying the
 set $A_n$ of {\em compositions} into at most $n$ parts satisfying
\[
\frac{\la_1}{1}\ge\frac{\la_2}{2}\ge \ldots \ge \frac{\la_n}{n}\ge 0.
\] 
We refered to these as
{\em anti-lecture hall compositions} and showed the following
with a bijective proof along the lines of  Yee's proof of (\ref{lhpuv})
in \cite{yee2}

\noindent
{\rm \bf  The Refined Anti-Lecture Hall Theorem \cite{CS}~:}
\begin{equation}
A_n(u,v,q) \triangleq \sum_{\lambda\in A_n} q^{|\lambda|}
u^{|\lfloor \la \rfloor|}v^{o(\lfloor \la \rfloor)}=
\frac{(-uvq;q)_n}{(u^2q^{2};q)_n},
\label{alhuv}
\end{equation}
where $\lfloor \la \rfloor=(\lfloor \la_1/1\rfloor,
\lfloor \la_2/2\rfloor,\ldots,
\lfloor \la_n/n\rfloor)$ and $o(\la )$ denotes the number
of odd parts of a composition $\la$.

Setting $u=v=1$ in (\ref{alhuv}) gives the following analog of (\ref{lhp})~:

\noindent
{\bf The  Anti-Lecture Hall Theorem \cite{CS}~:}
\begin{equation}
A_n(q)  \triangleq 
\sum_{\lambda \in A_n} q^{|\lambda|} = \frac{(-q;q)_n}{(q^2;q)_n}.
\label{alh}
\end{equation}
At the time of \cite{CS} we did not have a $q$-series proof of
(\ref{alhuv}) and we were not able to prove an odd/even
refinement of  (\ref{alh}) analogous to (\ref{lhpxy}),
although we had a conjecture as to its form. 

In this paper we extend all of these results.
 Our starting
point is a proof
of The Refined Lecture Hall Theorem from  \cite{BME3},
where
Bousquet-M\'elou
and Eriksson gave a two-step  proof of (\ref{lhpuv})
using basic $q$-series identities. 
They explained separately
the numerator and the denominator using elementary techniques,
deriving a recurrence to obtain the denominator and noting that the
recurrence could be solved using a special case of the
$q$-analog of the Chu-Vandermonde summation (\cite{MR58:27738}, 3.3.10).

In Section 2.1, we pursue their approach, but proceed directly to a recurrence
for $L_n(u,v,q)$  which can be solved in a straightforward way
using the identity
$q$-Chu Vandermonde II \cite{gara}~:
$$
\frac{a^n(c/a;q)_n}{(c;q)_n}=
\sum_{m=0}^n \frac{(a;q)_m(q^{-n};q)_m}
{(c;q)_m(q;q)_m}q^m.
$$
With this modified approach, we give, in Section 2.2,
a new proof of the
Refined Anti-Lecture Hall  Theorem.  We show that it can be obtained
using another $q$-Chu Vandermonde I \cite{gara}~:
$$
\frac{(c/a;q)_n}{(c;q)_n}=
\sum_{m=0}^n \left[ \begin{matrix} n \\ m
\end{matrix} \right]_q    \frac{(a;q)_m}{(c;q)_m}(-c/a)^mq^{m\choose 2},
$$
where $\left[ \begin{matrix} n \\ m
\end{matrix} \right]_q = (q^{n-m+1};q)_m/(q;q)_m$
is the classical Gaussian polynomial,
the generating function for partitions into $m$ nonnegative
parts of size at most $n-m$.
  
    From this point on, all of our results are new.  We conjectured them
thanks to the 
Maple implementation of 
the generating function developed in \cite{CS1}.
In Section 3,
we show how the $q$-series techniques of Section 2 can be extended to get 
new identities for the enumeration of {\em truncated objects}.
For $n\ge k$, let $L_{n,k}$ be the
 set  of partitions
 $\lambda = (\lambda_1, \lambda_2, \ldots, \lambda_k)$ into  $k$ nonnegative parts satisfying
\[
\frac{\la_1}{n}\ge\frac{\la_2}{n-1}\ge \ldots \ge \frac{\la_k}{n-k+1}\ge 0.
\] 
We refer to these as
{\em truncated lecture hall partitions}.
(Note that the case $n=k$ corresponds to the ordinary
lecture hall partitions.)
These objects were introduced by  Eriksen in \cite{ne}, where
he gave a recurrence for their generating function
in one variable, but no closed-form solution.

We show in Section 3.1 how to compute
the three-variable generating function~:
\[
L_{n,k}(u,v,q) \triangleq   \sum_{\lambda\in L_{n,k}} q^{|\lambda|}
u^{|\lceil \la \rceil|}v^{o(\lfloor \la \rfloor)}
\]
where $\lceil \la \rceil = 
(\lceil \la_1/n \rceil,
\lceil \la_2/(n-1) \rceil, \ldots,
\lceil \la_{k-1}/(n-k+2) \rceil, 
\lceil \la_k/(n-k+1) \rceil)$. 

\noindent
{\rm \bf The Refined Truncated Lecture Hall Theorem:}
\begin{equation}
L_{n,k}(u,v,q) =
\sum_{m=0}^k(uv)^mq^{m+1\choose 2}\left[\begin{matrix}n\\ m\end{matrix}\right]_q
\frac{(-(u/v)q^{n-m+1};q)_m}{(u^2q^{2n-m+1};q)_m}.
\label{tlhuv}
\end{equation}

In Section 3.2 we study  {\em truncated anti-lecture hall compositions},
defined for each $n \geq k-1$ as the set $A_{n,k}$ 
of compositions
$\lambda = (\lambda_1, \lambda_2, \ldots, \lambda_k)$ 
 into $k$ nonnegative parts satisfying
\[
\frac{\la_1}{n-k+1}\ge\frac{\la_2}{n-k+2}\ge \ldots \ge \frac{\la_k}{n}\ge 0.
\]
(When $n=k$, these are the ordinary anti-lecture hall compositions.)
If we let
$\lfloor \la \rfloor$ denote the partition
$(\lfloor \la_1/(n-k+1) \rfloor,
\lfloor \la_2/(n-k+2) \rfloor, \ldots,
\lfloor \la_{k-1}/(n-1) \rfloor,
\lfloor \la_n/n \rfloor)$,
we  get the following, for $n\ge k$~:

\noindent
{\bf The Refined Truncated Anti-Lecture Hall Theorem:}
\begin{equation}
A_{n,k}(u,v,q) \triangleq   \sum_{\lambda\in A_{n,k}} q^{|\lambda|}
u^{|\lfloor \la \rfloor|}v^{o(\lfloor \la \rfloor)}=
\left[\begin{matrix}n\\ k\end{matrix}\right]_q
\frac{(-uvq^{n-k+1};q)_k}{(u^2q^{2(n-k+1)};q)_k}.
\label{talhuv}
\end{equation}

We next consider odd/even refinements of the generating
functions for truncated objects, analogous to (\ref{lhpxy}) for lecture
hall partitions. In Section 4.1, we show the following for
truncated lecture hall partitions.

\noindent
{\bf The Odd/Even Truncated Lecture Hall Theorem:}
\begin{equation}
L_{n,k}(x,y)  \triangleq 
\sum_{\lambda\in L_{n,k}}x^{|\lambda_o|}y^{|\lambda_o|}=
\sum_{m=0}^k
\frac{
\left(x^{\lfloor m/2\rfloor+1}y^{\lfloor m/2\rfloor}\right)^{\lceil m/2\rceil}
\left[\begin{matrix} n-\lceil m/2\rceil\\ 
\lfloor m/2\rfloor\end{matrix}\right]_{xy}}
{(x;xy)_{\lceil m/2\rceil}(x^ny^{n-1};(xy)^{-1})_{\lfloor m/2\rfloor}}.
\label{tlhxy}
\end{equation}

Similarly, in Section 4.2 we find the odd/even generating
function  for truncated
anti-lecture hall compositions. Let $n\ge k-1$.

\noindent
{\bf The Odd/Even Truncated Anti-Lecture Hall Theorem:}
\begin{equation}
A_{n,k}(x,y) \triangleq 
\sum_{\lambda\in A_{n,k}}x^{|\lambda_o|}y^{|\lambda_o|}=
\frac{\left[\begin{matrix}n\\ \lfloor k/2\rfloor \end{matrix}\right]_{xy}}
{(x;xy)_{\lceil k/2\rceil }(x^{n-k+1}y^{n-k+2};xy)_{\lfloor k/2\rfloor}}.
\label{talhxy}
\end{equation}

In particular, setting $n=k$ in (\ref{talhxy})
  gives for the first time the
odd/even generating function for anti-lecture hall compositions~:

\begin{equation}
A_n(x,y)=A_{n,n}(x,y)=\frac{\left[\begin{matrix}n\\ \lfloor n/2\rfloor \end{matrix}\right]_{xy}}
{(x;xy)_{\lceil n/2\rceil}(xy^{2};xy)_{\lfloor n/2\rfloor}}.
\label{2alhc}
\end{equation}

One of the many interesting things about lecture hall partitions
is that the Lecture Hall Theorem gives them  a simple interpretation
in terms of partitions into odd parts:
{\em the number of partitions of $N$ in $L_n$
is equal to the number of partitions of $N$ into odd parts
less than $2n$.}  A similar interpretation of truncated lecture
hall partitions is not so evident from their generating function
in (\ref{tlhuv}) or in (\ref{tlhxy}). 
However, in Section 5 we show  the following correspondence between
truncated lecture hall partitions and partitions into odd parts 
with certain restrictions.

\noindent
{\bf Characterization of Truncated
Lecture Hall Partitions:}
\begin{quote}
%For $n \geq k \geq 1$, 
The number of truncated lecture hall partitions
of $N$ in $L_{n,k}$ is equal to the number of partitions of $N$ into 
odd parts less than $2n$, with the following constraint on the parts:
at most  
$\lfloor k/2\rfloor$
parts can be chosen from the set
$$
\{2\lceil k/2\rceil +1,2\lceil k/2\rceil +3,\ldots ,
2(n-\lfloor k/2\rfloor)-1\}.$$
\end{quote}

As $n \rightarrow \infty$, the set $L_n$ of lecture hall partitions
approaches the set of partitions into distinct parts.
Similarly, the right-hand side of (\ref{lhp}) approaches the set
of partitions into odd parts. In this sense, the Lecture Hall Theorem
is viewed as a {\em finitization} of Euler's Theorem, which states:
{\em The number of
partitions of an integer $N$ into distinct parts is equal to the
number of partitions of $N$ into odd parts.}
In Section 5.1, we show how our results on truncated lecture hall
partitions lead
to finizations of certain  refinements of Euler's Theorem
that are implied by Sylvester's bijection.

We show in Section 5.2 how  truncated anti-lecture hall theorems
can be viewed as finitizations of another well-known identity:
{\em The number of partitions of $N$ into $k$ parts is equal to the number
of partitions of $N$ with no part larger than $k$.}

\section{The Refined Lecture Hall Theorems}

We will make use of the $q$-multinomal coefficient, defined for
$n=n_0+n_1+ \cdots + n_t$ by
\[
\left[
\begin{matrix}
n\\
n_0,n_1,\ldots, n_t\end{matrix}\right]_q =
\frac{(q;q)_n}{(q;q)_{n_0}(q;q)_{n_1} \cdots (q;q)_{n_t}},
\] from which it follows that
\begin{equation}
\left[
\begin{matrix}
n\\
n_0,n_1\ldots n_t\end{matrix}\right]_q =
\left[
\begin{matrix}
n\\
n_0\end{matrix}\right]_q 
\left[
\begin{matrix}
n-m_0\\
n_1,\ldots, n_t\end{matrix}\right]_q. 
\label{multinom}
\end{equation}

\subsection{Lecture Hall Partitions}

We review the proof of the Refined Lecture Hall Theorem
(\ref{lhpuv}) of 
Bousquet-M{\'e}lou and Eriksson.
Most of the basic ideas of this proof come from \cite{BME3}.
Given a lecture hall partition $\lambda \in L_n$,
denote by $\lceil \lambda\rceil$ the sequence
$(\lceil \lambda_1/n\rceil,\lceil \lambda_2/(n-1)\rceil,\ldots , \lceil \lambda_n/1\rceil)$.
We can write  %$\lambda$ as
%$$((\mu_1,\ldots ,\mu_n),(r_1,\ldots ,r_n))$$ where
$\lambda_i=(n-i+1)\mu_i-r_i$, with $0\le r_i\le n-i$ for $1\le i\le n$.
Then $(\mu_1, \ldots, \mu_n)=\lceil \lambda \rceil$.
\begin{proposition}\cite{BME3} 
A partition $\lambda$ is in  $L_n$ if and only if 
\begin{enumerate}
\item $\mu_1\ge \mu_2\ge \ldots \ge \mu_n\ge 0$  and\\  
\item $r_i\le r_{i+1}$ whenever $\mu_i=\mu_{i+1}$.
\end{enumerate}
\label{cond1}
\end{proposition}
\noindent
The first condition  implies that $\lceil \lambda \rceil$  is a partition
into $n$ nonnegative parts.

Let $P_n$ be the set of partitions into $n$ nonnegative parts.
Let $\mu$ be a partition in $P_n$.
In \cite{BME3} the authors compute the generating
function of the lecture hall partitions $\lambda$
having $\lceil \lambda \rceil=\mu$.
Let this generating function
be $L_\mu(q)$.
\begin{proposition}{\cite{BME3}}
For $\mu \in P_n$,
\[
L_{\mu}(q) \triangleq 
\sum_{\begin{subarray}{c}\lambda\in L_n\\
\lceil \lambda\rceil =\mu\end{subarray}}q^{|\lambda|} =
q^{\sum_{i=1}^n (n-i+1)\mu_i}\left[
\begin{matrix}
n\\
m_0,m_1,\ldots, m_{\mu_1}\end{matrix}\right]_{1/q}.
\]
where $m_i$ is the multiplicity of the part $i$ in $\mu$.
\label{Lu1}
\end{proposition}
\noindent{\bf Proof.} We sketch here the main ideas of the proof. 
As $|\lambda|=\sum_{i=1}^n ((n-i+1)\mu_i-r_i)$, 
\[
L_\mu(q)=q^{\sum_{i=1}^n (n-i+1)\mu_i}\sum_{(r_1,\ldots, r_n)}q^{-\sum_{i=1}^n 
r_i}. 
\]
For $0 \leq i \leq \mu_1$,
let $$\ell_i=\left\{\begin{matrix} 0&{\rm if}\ i=0\\ 
\sum_{j=\mu_1-i+1}^{\mu_1} m_{j}& {\rm otherwise.}
\end{matrix}\right.$$ 
Then
the condition (2) in Proposition \ref{cond1} implies that 
$(r_{\ell_{i+1}},r_{\ell_{i+1}-1},\ldots ,r_{\ell_{i}+1})$ is a partition
into $\ell_{i+1}-\ell_{i}=m_{\mu_1-i}$  nonnegative parts and that these parts
are less than or equal to $n-\ell_{i+1}$.
Therefore their generating function is well-known to be a Gaussian polynomial~:
\begin{equation}
\sum_{(r_{\ell_{i+1}},\ldots ,r_{\ell_{i}+1})}
q^{(r_{\ell_{i+1}}+\ldots +r_{\ell_{i}+1})}
=\left[ \begin{matrix} n-\ell_{i+1}+m_{\mu_1-i} \\m_{\mu_1-i}
\end{matrix}\right]_{q}
=\left[ \begin{matrix} n-\ell_i\\ m_{\mu_1-i}\end{matrix}\right]_{q}.
\label{gaussgf}
\end{equation}
Hence, the result follows from the computation below which uses
(\ref{gaussgf}) for the second equality and repeated application of
 (\ref{multinom}) for the third.
\begin{eqnarray*}
\sum_{(r_1,\ldots ,r_n)}q^{-\sum_{i=1}^n r_i}
&=&\prod_{i=0}^{\mu_1} \sum_{(r_{\ell_{i+1}},\ldots ,r_{\ell_{i}+1})}
(1/q)^{(r_{\ell_{i+1}}+\ldots +r_{\ell_{i}+1})}\\
&=&\prod_{i=0}^{\mu_1}\left[ \begin{matrix} 
n-\ell_i\\ m_{\mu_1-i}\end{matrix}\right]_{1/q}\\
%&=&\prod_{i}\left[ \begin{matrix} n-\ell_i\\ m_{\mu_1-i}\end{matrix}\right]_{1/q}\\
%&=&\left[ \begin{matrix} n\\ 
%m_{\mu_1},m_{\mu_1-1},\ldots ,m_0\end{matrix}\right]_{1/q}\\
&=& \left[ \begin{matrix} n\\ m_{0},m_{1},\ldots ,m_{\mu_1}
\end{matrix}\right]_{1/q}.
\end{eqnarray*}
\hfill $\Box$

At this point, we take a different turn from \cite{BME3},
proceeding directly to enumeration of the partitions $\mu$ with
all parts positive.
Let $P_{n,m}$ be the set of partitions into   $n$ nonnegative parts,
$m$ of which are positive.
Given $\mu$ in $P_{m,n}$, we define $\tilde{\mu}$ in $P_m$  
by  $\tilde{\mu}_i=\mu_i-1$,  $1\le i\le m$.
Then
\begin{proposition} For $\mu \in P_{m,n}$,
\[
L_{\mu}(q)=q^{(n-m)|\tilde{\mu}|+{m+1\choose 2}}\left[
\begin{matrix}
n\\
m\end{matrix}\right]_{q}L_{\tilde{\mu}}(q).
\]
\label{Lu2}
\end{proposition}
\noindent{\bf Proof.}
Using  Proposition \ref{Lu1}, we know that if $m_i$ is the multiplicity of $i$ in $\mu$, then
\[
L_\mu(q)=q^{\sum_{i=1}^m (n-i+1)(\tilde{\mu}_i+1)}\left[ \begin{matrix} n\\ m_{0},m_{1},\ldots ,m_{\mu_1}
\end{matrix}\right]_{1/q},
\]
and that
\[
L_{\tilde{\mu}}(q)=q^{\sum_{i=1}^m (m-i+1)\tilde{\mu}_i}\left[ \begin{matrix} m\\ m_{1},\ldots ,m_{\mu_1}
\end{matrix}\right]_{1/q}.
\]
Since $m_0=n-m$, using (\ref{multinom}), we can write the first equation as
%As $\left[ \begin{matrix} n\\ m_{0},m_{1},\ldots ,m_{\mu_1}
%\end{matrix}\right]_{1/q}=
%\left[ \begin{matrix} n\\ m_0\end{matrix}\right]_{1/q}
%\left[ \begin{matrix} m\\ m_{1},\ldots ,m_{\mu_1}
%\end{matrix}\right]_{1/q}$, we get ~:
\[
L_\mu(q)=q^{(n-m)|\tilde{\mu}|+m(n+1)-m(m+1)/2}\left[ \begin{matrix} n\\ m\end{matrix}\right]_{1/q}q^{\sum_{i=1}^m (m-i+1)\tilde{\mu}_i}\left[ \begin{matrix} m\\ m_{1},\ldots ,m_{\mu_1}
\end{matrix}\right]_{1/q}.
\]
The result follows from the second equation and the identity
\begin{equation*}
\left[ \begin{matrix} n\\ m\end{matrix}\right]_{1/q} =
q^{-m(n-m)}\left[ \begin{matrix} n\\ m\end{matrix}\right]_{q}.
%\label{1overq}
\end{equation*}
\hfill $\Box$\\

An easy consequence of Proposition \ref{Lu2} is that
\begin{equation}
u^{|\mu|}v^{o(\mu)}L_{\mu}(q)=(uv)^mu^{|\tilde{\mu}|}v^{-o(\tilde{\mu})}q^{(n-m)|\tilde{\mu}|+{m+1\choose 2}}\left[
\begin{matrix}
n\\
m\end{matrix}\right]_{q}L_{\tilde{\mu}}(q),
\label{need}
\end{equation}
as $|\mu|=|\tilde{\mu}|+m$ and $o(\mu)=m-o(\tilde{\mu})$.

In order to prove the Refined Lecture Hall Theorem
(\ref{lhpuv}), the generating function we are looking for,
$L_n(u,v,q) \triangleq  
\sum_{\lambda\in L_n }u^{|\lceil \lambda\rceil|}
v^{o({\lceil\lambda\rceil})}q^{|\lambda|}$,  can be rewritten  as~:
\begin{equation}
L_n(u,v,q)=\sum_{\mu\in P_n}u^{|\mu|}v^{o(\mu)}L_\mu(q)=\sum_{m=0}^n \sum_{\mu\in P_{n,m}} u^{|\mu|}v^{o(\mu)}L_{\mu}(q).
\label{Lndef}
\end{equation}
%where $P_{n,m}$ is the set of partitions into $n$ nonnegative parts, $m$
%of which are positive.
Then we can prove the following recurrence for $L_n(u,v,q)$.
\begin{proposition}  $L_0(u,v,q) = 1$ and for $n>0$,
\[
L_n(u,v,q)=\sum_{m=0}^n \left[
\begin{matrix}
n\\
m\end{matrix}\right]_q (uv)^m q^{m+1\choose 2}L_m(uq^{n-m},1/v,q).
\]
\label{Lu3}
\end{proposition}
\noindent{\bf Proof.} 
\begin{eqnarray*}
L_n(u,v,q)%&=&\sum_{\mu\in P_n}u^{|\mu|}v^{o(\mu)}L_\mu\\
&=&\sum_{m=0}^n\sum_{\mu\in P_{n,m} 
}u^{|\mu|}v^{o(\mu)}L_{\mu}(q)\\
&=&\sum_{m=0}^n \sum_{\tilde{\mu}\in P_m}(uv)^mu^{|\tilde{\mu}|}v^{-o(\tilde{\mu})}
q^{(n-m)|\tilde{\mu}|+{m+1\choose 2}}\left[\begin{matrix}
n\\m\end{matrix}\right]_{q}L_{\tilde{\mu}}(q) \ \ \ \ 
({\rm using}\ (\ref{need}))\\
&=&\sum_{m=0}^n \left[
\begin{matrix}
n\\
m\end{matrix}\right]_{q}(uv)^mq^{{m+1\choose 2}}\sum_{\tilde{\mu}\in P_m}
(uq^{n-m})^{|\tilde{\mu}|}v^{-o(\tilde{\mu})}L_{\tilde{\mu}}(q)
\end{eqnarray*}
If we now apply the first equality of (\ref{Lndef}) to
$L_m(uq^{n-m},1/v,q)$, we get the last sum in the last line
above and the proposition is proved.
%
%&=&\sum_{m=0}^n \left[
%\begin{matrix}
%n\\
%m\end{matrix}\right]_q (uv)^m q^{m+1\choose 2}L_m(uq^{n-m},1/v,q).
%\end{eqnarray*}
\hfill $\Box$

\noindent{\bf Proof of The Refined Lecture Hall Theorem}
(\ref{lhpuv}):
We solve the recurrence of Proposition \ref{Lu3}
using the identity
$q$-Chu Vandermonde II:
\begin{equation}
\frac{a^n(c/a;q)_n}{(c;q)_n}=
\sum_{m=0}^n \frac{(a;q)_m(q^{-n};q)_m}
{(c;q)_m(q;q)_m}q^m.
\label{qchuII}
\end{equation}
If we set $a=-vq^{-n}/u$, $c=q^{-2n}/u^2$ in (\ref{qchuII}), we get that
$$
(-vq^{-n}/u)^n\frac{(-q^{-n}/uv;q)_n}{(q^{-2n}/u^2;q)_n}=
\sum_{m=0}^n \frac{(-vq^{-n}/u;q)_m(q^{-n};q)_m}
{(q^{-2n}/u^2;q)_m(q;q)_m}q^m.
$$
Now on every factor above of the form $(bq^{-t};q)_m$ we use the
identity
\[
(bq^{-t};q)_m=(-bq^{-t+(m-1)/2})^m(q^{t-m+1}/b;q)_m
\]
and get
$$
\frac{(-uvq;q)_n}{(u^2q^{n+1};q)_n}=
\sum_{m=0}^n
(uv)^m q^{m+1\choose 2}\left[
\begin{matrix}
n\\
m\end{matrix}\right]_{q} \frac{(-(u/v)q^{n-m+1};q)_m}{(u^2q^{2n-m+1};q)_m}.
$$
This shows that
$$L_n(u,v,q)=\frac{(-uvq;q)_n}{(u^2q^{n+1};q)_n}$$
is the solution to the recurrence of 
Proposition \ref{Lu3}.
\hfill $\Box$

\subsection{Anti-Lecture Hall Compositions}

A $q$-series proof of the Refined Anti-Lecture Hall Theorem 
(\ref{alhuv}) will follow the same
approach as in the previous subsection.
Given an anti-lecture hall composition $\lambda \in A_n$,
define the floor of   $\lambda$, as 
$\lfloor \lambda\rfloor=(\lfloor \lambda_1/1\rfloor,
\lfloor \lambda_2/2\rfloor,\ldots ,\lfloor \lambda_n/n\rfloor).$ 
Then write %$\lambda$ ad
%$$((\mu_1,\ldots ,\mu_n),(r_1,\ldots ,r_n))$$ where
$\lambda_i=i\mu_i+r_i$, with $0\le r_i\le i-1$ for $1\le i\le n$.
Note that $(\mu_1, \ldots, \mu_n)=\lfloor \lambda \rfloor$
and that $\lambda\in A_n$ if and only if \\
(i) $\mu_1\ge \mu_2\ge \ldots \ge \mu_n\ge 0$ and \\ 
(ii) $r_i\ge r_{i+1}$ whenever $\mu_i=\mu_{i+1}$.

The condition (i) implies that $\lfloor \lambda\rfloor$ 
is a partition
in $P_n$.
We fix $\mu \in P_n$ and
compute the generating
function $A_{\mu}$ of the anti-lecture hall compositions $\lambda$
having  $\lfloor \lambda \rfloor=\mu$.
%$A_\mu=
%\sum_{
%\begin{subarray}{c} \lambda\in A_n \\
%\lfloor \lambda\rfloor =\mu
%\end{subarray}}q^{|\lambda|}$.
\begin{proposition} For $\mu \in P_n$,
\[
A_{\mu}(q)  \triangleq 
\sum_{
\begin{subarray}{c} \lambda\in A_n \\
\lfloor \lambda\rfloor =\mu
\end{subarray}}q^{|\lambda|}
= q^{\sum_{i=1}^n i\mu_i}\left[
\begin{matrix}
n\\
m_0,m_1,\ldots, m_{\mu_1}\end{matrix}\right]_{q}.
\]
where $m_i$ is the multiplicity of the part $i$ in $\mu$.
\label{Ans}
\end{proposition}
\noindent{\bf Proof.}
As $\lambda_i=i\mu_i+r_i$
\[
A_{\mu}(q)=q^{\sum_{i=1}^n i\mu_i}\sum_{(r_1,\ldots ,r_n)}q^{\sum_{i=1}^n r_i}.
\]  For $0 \leq i \leq \mu_1+1$,
let $\ell_i= 
n-\sum_{j=0}^{i-1} m_{j}$.
Then
the condition (ii) implies that 
$(r_{\ell_{i+1}+1},r_{\ell_{i+1}+2},\ldots ,r_{\ell{i}})$ is a partition
into $\ell_{i}-\ell_{i+1}=m_{i}$  nonnegative parts and that these parts
are less than  or equal to $\ell_{i+1}$.
%Therefore their generating function is well known to be
%$$
%\sum_{(r_{\ell_{i+1}},\ldots ,r_{\ell_{i}+1})}
%q^{(r_{\ell_{i+1}}+\ldots +r_{\ell_{i}+1})}
%=\left[ \begin{matrix} \ell_{i}\\ m_{i}\end{matrix}\right]_{q}
%$$
Therefore, since $\ell_0=n$,
$$
\sum_{(r_1,\ldots ,r_n)}q^{\sum_{i=1}^n r_i}
=\prod_{i=0}^{\mu_1} 
\left[ \begin{matrix} \ell_{i+1}+m_i\\ m_i\end{matrix}\right]_{q}\\
= \left[ \begin{matrix} n\\ m_{0},m_{1},\ldots ,m_{\mu_1}
\end{matrix}\right]_{q}.  
$$
\hfill $\Box$

As before
if $\mu \in P_{n,m}$  then we can get
$\tilde{\mu}$ in $P_m$ by  $\tilde{\mu}_i=\mu_i-1$,  $1\le i\le m$.
\begin{proposition}  For $\mu \in P_{n,m}$,
\[
A_{\mu}(q)=q^{{m+1\choose 2}}\left[
\begin{matrix}
n\\
m\end{matrix}\right]_{q}A_{\tilde{\mu}}(q).
\]
\label{Au2}
\end{proposition}
\noindent{\bf Proof.} From  Proposition \ref{Ans} we know that
\[
A_\mu(q)=q^{\sum_{i=1}^m i(\tilde{\mu}_i+1)}\left[ \begin{matrix} n\\ m_{0},m_{1},\ldots ,m_{\mu_1}
\end{matrix}\right]_{q},
\]
and that
\[
A_{\tilde{\mu}}(q)=q^{\sum_{i=1}^m i\tilde{\mu}_i}\left[ \begin{matrix} m\\ m_{1},\ldots ,m_{\mu_1}
\end{matrix}\right]_{q}.
\]
Since $m=n-m_0$, we get the result. \hfill $\Box$\\

If $\mu\in P_{n,m}$, then
 $|\mu|=|\tilde{\mu}|+m$ and $o(\mu)=m-o(\tilde{\mu})$, so
we can conclude from Proposition \ref{Au2}
that~: 
\begin{equation}
u^{|\mu|}v^{o(\mu)}A_{\mu}(q)=(uv)^mu^{|\tilde{\mu}|}v^{-o(\tilde{\mu})}q^{{m+1\choose 2}}\left[
\begin{matrix}
n\\
m\end{matrix}\right]_{q}A_{\tilde{\mu}}(q).
\label{Au3}
\end{equation}

In order to prove the Refined Anti-Lecture Hall Theorem
(\ref{alhuv}), the generating function we are looking for,
$A_n(u,v,q) \triangleq 
\sum_{\lambda\in A_n }u^{|\lfloor \lambda\rfloor|}
v^{o({\lfloor\lambda\rfloor})}q^{|\lambda|}$
can be rewritten
as
\begin{equation}
A_n(u,v,q) =
\sum_{\mu \in P_n}u^{|\mu|}
v^{o(u)}A_{\mu}(q) =
 \sum_{m=0}^n\sum_{\mu\in P_{n,m}}u^{|\mu|}v^{o(\mu)}A_\mu(q).
\label{An1}
\end{equation}
Then we can prove the following recurrence for $A_n(u,v,q)$.
\begin{proposition}
 $A_0(u,v,q) = 1$ and for $n>0$
\[
A_n(u,v,q)=\sum_{m=0}^n \left[
\begin{matrix}
n\\
m\end{matrix}\right]_q (uv)^m q^{m+1\choose 2}A_m(u,1/v,q).
\]
\label{prop2}
\end{proposition}
\noindent{\bf Proof.} 
\begin{eqnarray*}
 A_n(u,v,q)%&=&\sum_{\mu\in P_n}u^{|\mu|}v^{o(\mu)}A_\mu(q)\\
&=&
\sum_{m=0}^n\sum_{\mu\in P_{n,m}}u^{|\mu|}v^{o(\mu)}A_{\mu}(q)\\
&=&\sum_{m=0}^n \sum_{\tilde{\mu}\in P_m}(uv)^mu^{|\tilde{\mu}|}v^{-o(\tilde{\mu})}
q^{{m+1\choose 2}}\left[\begin{matrix}
n\\m\end{matrix}\right]_{q}A_{\tilde{\mu}}(q) 
\ \ \ \
({\rm using}\ (\ref{Au3}))\\
%&=&\sum_{m=0}^n \left[
%\begin{matrix}
%n\\
%m\end{matrix}\right]_{q}(uv)^mq^{{m+1\choose 2}}\sum_{\tilde{\mu}\in P_m}
%(u)^{|\tilde{\mu}|}v^{-o(\tilde{\mu})}A_{\tilde{\mu}}\\
&=&\sum_{m=0}^n \left[
\begin{matrix}
n\\
m\end{matrix}\right]_q (uv)^m q^{m+1\choose 2}A_m(u ,1/v,q),
\end{eqnarray*}
where the last step follows from
first equality of (\ref{An1}).
\hfill $\Box$

\noindent
{\bf Proof of the Refined Anti-Lecture Hall Theorem} (\ref{alhuv}).
We solve the recurrence of Proposition \ref{prop2} using the
identity
$q$-Chu Vandermonde I:
\begin{equation}
\frac{(c/a;q)_n}{(c;q)_n}=
\sum_{m=0}^n \left[ \begin{matrix} n \\ m
\end{matrix} \right]_q    \frac{(a;q)_m}{(c;q)_m}(-c/a)^mq^{m\choose 2}.
\label{qchuI}
\end{equation}
Setting $a= -uq/v$ and  $c = u^2q^2$ in (\ref{qchuI}) 
gives immediately
$$
\frac{(-uvq;q)_n}{(u^2q^{2};q)_n}=
\sum_{m=0}^n
\left[ \begin{matrix} n\\ m\end{matrix}\right]_{q}
(uv)^m q^{m+1\choose 2}
 \frac{(-uq/v;q)_m}{(u^2q^{2};q)_m}.
$$
With the recurrence of Proposition \ref{prop2}, we conclude that~:
$$A_n(u,v,q)=\frac{(-uvq;q)_n}{(u^2q^{2};q)_n}.$$
 \hfill $\Box$

\section{Refined Lecture Hall Theorems for Truncated Objects}

In this section we apply the techniques used in Section 2
to derive the new refined truncated lecture hall theorems.

\subsection{Truncated Lecture Hall Partitions}
Recall that a truncated lecture hall partition in $L_{n,k}$
is a sequence $\lambda=(\lambda_1,\ldots ,\lambda_k)$ such that
\[
\frac{\lambda_1}{n}\ge \frac{\lambda_2}{n-1}\ge
\ldots \ge \frac{\lambda_k}{n-k+1}\ge 0.
\] 
Given a lecture hall partition $\lambda \in L_{n,k}$,
 we write  
$\lambda_i=(n-i+1)\mu_i-r_i$, with $0\le r_i\le n-i$ for $1\le i\le k$.
Let $\lceil \lambda \rceil=(\lceil \lambda_1/n \rceil, \ldots ,
\lceil \lambda_k/(n-k+1) \rceil)$. Then
$\mu=\lceil \lambda \rceil$. Note that
$\lambda$ has $k$ positive parts if and only if  $\mu$
has $k$ positive parts.

\begin{theorem}
{\rm  (The Refined Truncated Lecture Hall Theorem)}
\begin{equation}
L_{n,k}(u,v,q)=
\sum_{m=0}^k(uv)^mq^{m+1\choose 2}\left[\begin{matrix}n\\ m\end{matrix}\right]_q\frac{(-(u/v)q^{n-m+1};q)_m}{(u^2q^{2n-m+1};q)_m}.
\label{rtlh}
\end{equation}
%with $\left[\begin{matrix}k\\ m\end{matrix}\right]_q
%=(q^{k-m+1};q)_m/(q;q)_m$.
\end{theorem}

\noindent
{\bf Proof.}
Let $\bar{L}_{n,k}$ be the set of lecture hall partitions
in $L_n$ with $k$ positive parts.
Then
\[
\bar{L}_{n,k}(u,v,q) \triangleq 
\sum_{\lambda \in \bar{L}_{n,k}}
u^{|\lceil \lambda \rceil |}
v^{o(\lceil \lambda \rceil)} 
q^{|\lambda|}
= \sum_{\mu \in P_{n,k}} u^{|\mu|}v^{o(\mu)}L_{\mu}.
\]
It follows from the proof of Proposition \ref{Lu3} (ignoring the outer
sum there) that
\[
\bar{L}_{n,k}(u,v,q)=
(uv)^kq^{k+1\choose 2}\left[\begin{matrix}n\\ k\end{matrix}\right]_q
L_k(uq^{n-k},1/v,q)
\]
and then applying the Refined Lecture Hall Theorem (\ref{lhpuv})
to $L_k$ gives
\begin{equation}
\bar{L}_{n,k}(u,v,q)=
(uv)^kq^{k+1\choose 2}\left[\begin{matrix}n\\ k\end{matrix}\right]_q
\frac{(-(u/v)q^{n-k+1};q)_k}{(u^2q^{2n-k+1};q)_k}.
\label{rtlh=}
\end{equation}
The result follows since
 $L_{n,k}(u,v,q)
= \sum_{m=0}^k \bar{L}_{n,m}(u,v,q)$. 
\hfill $\Box$

\subsection{Truncated Anti-Lecture Hall Compositions}

Let $n\ge k$.
A truncated anti-lecture hall composition in $A_{n,k}$
is a sequence $\lambda=(\lambda_1,\ldots ,\lambda_k)$ such that
\[
\frac{\lambda_1}{n-k+1}\ge \frac{\lambda_2}{n-k+2}\ge 
\ldots \ge \frac{\lambda_k}{n}\ge 0.
\] 
We write
$\lambda_i=(n-k+i)\mu_i+r_i$, with $0\le r_i\le n-k+i-1$ for $1\le i\le k$
and define $\lfloor \lambda \rfloor=(\mu_1, \ldots, \mu_k)$.
As before,  $\lambda\in A_{n,k}$ if and only if \\
(i) $\mu_1\ge \mu_2\ge \ldots \ge \mu_k\ge 0$ and \\ 
(ii) $r_i\ge r_{i+1}$ whenever $\mu_i=\mu_{i+1}$.

  For  $\mu \in P_k$, we
compute the generating
function $A_{\mu,k}(q)$ of those $\lambda \in A_{n,k}$
having  $\lfloor \lambda \rfloor=\mu$. 
%$A_{\mu,n}=\sum_{\begin{subarray}{c}\lambda\in A_{n,k}\\ 
%\lfloor \lambda\rfloor =\mu\end{subarray}}q^{|\lambda|}$.
\begin{proposition}
\[
A_{\mu,k}(q)\triangleq 
\sum_{\begin{subarray}{c}\lambda\in A_{n,k}\\
\lfloor \lambda\rfloor =\mu\end{subarray}}q^{|\lambda|} =
\left[\begin{matrix}n\\ k\end{matrix}\right]_q q^{(n-k)|\mu|}A_{\mu}(q),
\]
where $A_{\mu}(q)$ is the generating  of $\lambda \in A_{k}$
having  $\lfloor \lambda \rfloor=\mu$. 
\label{Aun}
\end{proposition}
\noindent{\bf Proof.}
\[
A_{\mu,k}(q)=q^{\sum_{i=1}^k (n-k+i)\mu_i}\sum_{(r_1,\ldots ,r_k)}q^{\sum_{i=1}^k r_i}.
\] For $0 \leq i \leq \mu_1 + 1$,
let $\ell_i= 
n-\sum_{j=0}^{i-1} m_{j}$. 
Then
the condition (ii) implies that 
$(r_{\ell_{i+1}+1},r_{\ell_{i+1}+2},\ldots ,r_{\ell{i}})$ is a partition
into $m_{i}$  nonnegative parts and that these parts
are less than  or equal to $\ell_{i+1}$.
%Therefore their generating function is well known to be
%$$
%\sum_{(r_{\ell_{i+1}},\ldots ,r_{\ell_{i}+1})}
%q^{(r_{\ell_{i+1}}+\ldots +r_{\ell_{i}+1})}
%=\left[ \begin{matrix} \ell_{i}\\ m_{i}\end{matrix}\right]_{q}
%$$
Therefore, since $\ell_i-\ell_{i+1}=m_i$ and $\ell_{\mu_1+1}=n-k$,
\begin{equation*}
\sum_{(r_1,\ldots ,r_k)}q^{\sum_{i=1}^k r_i}
= \prod_{i=0}^{\mu_1} \left[ \begin{matrix} \ell_{i+1} + m_i\\ m_i
\end{matrix}\right]_{q} 
= \left[ \begin{matrix} n\\ m_{0},m_{1},\ldots ,m_{\mu_1},n-k
\end{matrix}\right]_{q}=
\left[ \begin{matrix} n\\k
\end{matrix}\right]_{q}
\left[ \begin{matrix} k\\ m_{0},m_{1},\ldots ,m_{\mu_1}
\end{matrix}\right]_{q}.
\end{equation*}
By Proposition \ref{Ans},
\[
A_{\mu}(q)=q^{\sum_{i=1}^k i\mu_i}\left[\begin{matrix} k\\
m_0,m_1, \ldots, m_{\mu_1}\end{matrix}\right]_q,
\]
%Now as $|\lambda|= (n-k)|\mu|+
%\sum_{i=1}^k i\mu_i+r_i$, 
and the result follows. \hfill $\Box$\\

\begin{theorem}
{\rm  (The Refined Truncated Anti-Lecture Hall Theorem)}.
\begin{equation}
A_{n,k}(u,v,q) \triangleq \sum_{\lambda\in A_{n,k}} q^{|\lambda|}
u^{|\lfloor \la \rfloor|}v^{o(\lfloor \la \rfloor)}=
\left[\begin{matrix}n\\ k\end{matrix}\right]_q
\frac{(-uvq^{n-k+1};q)_k}{(u^2q^{2(n-k+1)};q)_k}.
\label{rtalh}
\end{equation}
\label{rtalht}
\end{theorem}

\noindent{\bf Proof.} Applying first Proposition \ref{Aun}, then the
definition of $A_{k}(u,v,q)$,
and finally the Refined Anti-Lecture Hall
Theorem \ref{alhuv}, we get
\begin{eqnarray*}
A_{n,k}(u,v,q)
&=&\sum_{\mu\in P_k}
u^{|\mu|}v^{o(\mu)} 
\left[\begin{matrix}n\\ k\end{matrix}\right]_q q^{(n-k)|\mu|}A_{\mu}(q)\\
&=&\left[\begin{matrix}n\\ k\end{matrix}\right]_qA_k(uq^{n-k},v,q)\\
&=&\left[\begin{matrix}n\\ k\end{matrix}\right]_q\frac{(-uvq^{n-k+1};q)_k}
{(u^2q^{2(n-k+1)};q)_k}.
\end{eqnarray*}
\hfill $\Box$

\section{Odd/Even Generating Functions for Truncated Objects}

In this section we adapt the technique introduced in \cite{BME1}
to get the two-variable generating functions of the truncated 
objects. In \cite{BME1},
Bousquet-M{\'e}lou and Eriksson introduce a bijection
 BME: $L_{n-1}\times \mathbb{N} \rightarrow L_n$ that we recall here. 
For $\lambda\in L_{n-1}$ and $s\in \mathbb{N}$,
BME$(\lambda,s)=\mu$, where

\begin{eqnarray*}
\mu_1&\leftarrow &\left\lceil \frac{n\lambda_1}{n-1}\right\rceil+s \\
\mu_{2\ell} & \leftarrow &\lambda_{2\ell-1}, \ \ \ \ 1\le \ell\le  n/2;\\
\mu_{2\ell+1}&\leftarrow &\left\lceil \frac{(n-2\ell)\lambda_{2\ell+1}}{n-2\ell-1}\right\rceil+
\left\lfloor 
\frac{(n-2\ell)\lambda_{2\ell-1}}{n-2\ell+1}\right\rfloor-\lambda_{2\ell},
\ \ \ \ 1\le \ell\le  (n-1)/2. \\
\end{eqnarray*}

It is proved in  \cite{BME1} that $\mu\in {L}_{n}$, that
$|\mu_e|=|\lambda_o|$ and $|\mu_o|=2|\lambda_o|-|\lambda_e|+s$ and that BME
is a bijection.
This implies that~: 
\begin{equation*}
{L}_{n}(x,y)  \triangleq 
\sum_{\mu \in L_n} x^{|\mu_o|}y^{|\mu_e|}=
\sum_{\lambda \in L_{n-1}} 
    \sum_{s=0}^{\infty} x^{2|\lambda_o|-|\lambda_e|+s}y^{|\lambda_o|}=
\frac{1}{1-x}\sum_{\lambda \in L_{n-1}} 
     (x^2y)^{|\lambda_o|}(1/x)^{|\lambda_e|}=
\frac{{L}_{n-1}(x^2y,x^{-1})}{1-x},
\end{equation*}
giving the recurrence
\begin{equation}
{L}_{n}(x,y) =
\frac{{L}_{n-1}(x^2y,x^{-1})}{1-x}.
\label{Lnrec}
\end{equation}
As $L_0(x,y)=1$, this gives $L_n=1/(x;xy)_n$, the
Odd/Even Lecture Hall Theorem (\ref{lhpxy}).

\subsection{Truncated Lecture Hall Partitions}

Recall that $\bar{L}_{n,k}$ for $k\le n$ is the set
of partitions in $L_{n,k}$ with $k$ positive parts.
Let
\[
L_{n,k}(x,y) \triangleq \sum_{\lambda\in L_{n,k}}x^{|\lambda_o|}y^{|\lambda_e|}; \ \
\bar{L}_{n,k}(x,y) \triangleq \sum_{\lambda\in \bar{L}_{n,k}}x^{|\lambda_o|}y^{|\lambda_e|}.
\]
Note that $\bar{L}_{n,k}(x,y)=
{L}_{n,k}(x,y)-{L}_{n,k-1}(x,y)$.

For $n \geq k \geq 1$,
define a variation on the function BME, 
\[
{\rm BME}_{n,k}: \bar{L}_{n-1,k-1}\times \mathbb{N} \rightarrow \bar{L}_{n,k},
\]
by ${\rm BME}_{n,k}(\lambda,s) = (\mu_1, \mu_2, \ldots , \mu_k)$, where
%The new map  sends a partition $\lambda$ in 
%$\bar{L}_{n-1,2k}$ and an integer $s$ 
%into a partition $\mu$ in $\bar{L}_{n,2k+1}$~:
\begin{eqnarray*}
\mu_1&\leftarrow &\left\lceil \frac{n\lambda_1}{n-1}\right\rceil+s \\
\mu_{2\ell} & \leftarrow &\lambda_{2\ell-1}, \ \ 1\le \ell\le k/2;\\
\mu_{2\ell+1}&\leftarrow &\left\lceil \frac{(n-2\ell)\lambda_{2\ell+1}}{n-2\ell-1}\right\rceil+
\left\lfloor \frac{(n-2\ell)\lambda_{2\ell-1}}{n-2\ell+1}\right\rfloor
-\lambda_{2\ell},
\ \ 1\le \ell\le (k-1)/2; \\
\mbox{and, if $k$ is odd with $k=2t+1$,}\\
\mu_{2t+1}&\leftarrow &
\left\lfloor \frac{(n-2t)\lambda_{2t-1}}{n-2t+1}\right\rfloor-\lambda_{2t}+1.\\
\end{eqnarray*}
As was true for the function BME, 
${\rm BME}_{n,k}$ is one-to-one.
It is straightforward to check when $k$ is odd  that
the image of $ \bar{L}_{n-1,k-1}\times \mathbb{N}$ under
${\rm BME}_{n,k}$ is  $\bar{L}_{n,k}$, and that
$|\mu_e|=|\lambda_o|$ and $|\mu_o|=2|\lambda_o|-|\lambda_e|+s+1$.
This gives a recurrence for  $\bar{L}_{n,k}(x,y)$ when $k$ is odd.
\begin{proposition}
\begin{equation}
\bar{L}_{n,2k+1}(x,y) = 
\frac{x}{1-x}\bar{L}_{n-1,2k}(x^2y,x^{-1}).
\end{equation}
\label{map1}
\end{proposition}
Getting a recurrence for the even case will be harder.

For $i > 0$,
let $\bar{L}_{n,k,i}$ be the set of partitions in
$\bar{L}_{n,k}$ whose $k^{th}$ part is equal to $i$.
Let $\bar{L}_{n,k,0}=\bar{L}_{n,k-1}$.
We can check that when $k$ is even, ${\rm BME}_{n,k}$ gives a bijection
between  $\bar{L}_{n-1,k-1,i}\times \mathbb{N} $ and 
$\bar{L}_{n,k,i}$.  
%This other map will send a partition $\lambda$ in 
%$\bar{L}_{n-1,2k-1,i}$ and an integer $s$ 
%into a partition $\mu$ in $\bar{L}_{n,2k,i}$~:
%
%\begin{eqnarray*}
%\mu_1&\leftarrow &\lceil \frac{n\lambda_1}{n-1}\rceil+s \\
%%\mu_{2\ell} & \leftarrow &\lambda_{2\ell-1}; \ \ 1\le \ell\le 2k\\
%\mu_{2\ell+1}&\leftarrow &\lceil 
%\frac{(n-2\ell+1)\lambda_{2\ell+1}}{n-2\ell}\rceil+
%\lfloor \frac{(n-2\ell-1)\lambda_{2\ell-1}}{n-2\ell}\rfloor-\lambda_{2\ell};
%\ \ 1\le \ell\le k-1 \\
%\end{eqnarray*}
%
%It is easy to check that $\mu\in \bar{L}_{n,2k,i}$, that
Furthermore, when  $\lambda \in \bar{L}_{n-1,2k-1,i}$
and $\mu = {\rm BME}_{n,2k}(\lambda,s)$, then
$|\mu_e|=|\lambda_o|$ and 
$$|\mu_o|=
2|\lambda_o|-|\lambda_e|+s-\left\lfloor{\frac{(n-2k)i}{n-2k+1}}\right\rfloor
=2|\lambda_o|-|\lambda_e|+s+l-i,$$
with  $l=\lfloor i/(n-2k+1)\rfloor$. 
This implies that~: 
\begin{proposition}
For $i>0$,
\begin{equation}
\bar{L}_{n,2k,i}(x,y)=
\frac{x^{l-i}}{1-x}\bar{L}_{n-1,2k-1,i}(x^2y,1/x),
\end{equation}
with $l=\lfloor i/(n-2k+1)\rfloor$.
\label{bmeb}
\end{proposition}

Now we will decompose the set  $\bar{L}_{n,2k+1,i}$.
In what follows, in order to compress the notation, we will sometimes
write a function $f(x,y)$ as $f$, when the arguments are $(x,y)$.

\begin{proposition}
Let $i=l(n-2k)+r$, where
$1\le r \le n-2k$.

If $r \geq 2$, then 
\[
\bar{L}_{n,2k+1,i}=xL_{n,2k+1,i-1}-\frac{x}{1-x}\bar{L}_{n-1,2k-1,i+l}(x^2y,1/x).
\]
Otherwise, $r=1$ and
\[
\bar{L}_{n,2k+1,i}=xL_{n,2k+1,i-1}-
\frac{x}{1-x}\bar{L}_{n-1,2k-1,i+l-1}(x^2y,1/x)
-\frac{x}{1-x}\bar{L}_{n-1,2k-1,i+l}(x^2y,1/x).
\]
\label{decompose}
\end{proposition}
\noindent{\bf Proof.} For $i=l(n-2k)+r$,
\begin{equation}
\bar{L}_{n,2k+1,i}=\left\{
\begin{array}{ll}
x\bar{L}_{n,2k+1,i-1}-x^i\bar{L}_{n,2k,{i+l}} &{\rm if}\ r\ge 2 \ {\rm or}\ i=1\\
x\bar{L}_{n,2k+1,i-1}-x^i\bar{L}_{n,2k,i+l-1}-x^i\bar{L}_{n,2k,i+l}&{\rm if}\ r=1 \ {\rm and} \ i \not = 1\\
\end{array}\right.
\label{irecurr}
\end{equation}
which can be seen as follows.
Adding one to the $(2k+1)^{th}$ part of 
a partition in $\bar{L}_{n,2k+1,i-1}$ gives
a partition in $\bar{L}_{n,2k+1,i}$,
except in the following two cases (i) if the $2k^{th}$ part was equal to
$i+l$ and (ii)  if $r=1$ and the $2k^{th}$ part was equal to
$i+l-1$.

Now using Proposition \ref{bmeb},
for $i=l(n-2k)+r$ with $1\le r\le n-2k$,
\[
\bar{L}_{n,2k,i+l}=\frac{x^{1-i}}{1-x}\bar{L}_{n-1,2k-1,i+l}(x^2y,1/x),
\]
and for  $i=l(n-2k)+1$,
\[
\bar{L}_{n,2k,i+l-1}=\frac{x^{1-i}}{1-x}\bar{L}_{n-1,2k-1,i+l-1}(x^2y,1/x).
\]

Combining this with the recurrence (\ref{irecurr}),
we get the result. \hfill $\Box$

%Applying the recurrence, for $i=l(n-2k)+r$ and $2\le r\le n-2k$,
%LLL  n-1 should be n
%\[
%\bar{L}_{n,2k+1,i}=x\bar{L}_{n,2k+1,i-1}-x^i\bar{L}_{n-1,2k,i+l}.
%\]
%\begin{equation}
%\bar{L}_{n,2k+1,i}=x\bar{L}_{n,2k+1,i-1}-x^i\bar{L}_{n,2k,i+l}.
%\end{equation}
%RRR
%and using BME
%we get for $2\le r\le n-2k$,
%\[
%\bar{L}_{n,2k+1,i}=x\bar{L}_{n,2k+1,i-1}-\frac{x}{1-x}\bar{L}_{n-1,2k-1,i+l}(x^2y,1/x).
%\]

%Applying the recurrence, for $i=l(n-2k)+1$,
%LLL Both occurrences of n-1 should be n
%\[
%\bar{L}_{n,2k+1,i}=
%x\bar{L}_{n,2k+1,i-1}-x^i\bar{L}_{n-1,2k,i+l-1}-x^i\bar{L}_{n-1,2k,i+l}.
%\]
%\begin{equation}
%\bar{L}_{n,2k+1,i}=x\bar{L}_{n,2k+1,i-1}-x^i\bar{L}_{n,2k,i+l-1}-x^i\bar{L}_{n,2k,i+l}.
%\end{equation}
%ù%RRR
%we get
%\[
%\bar{L}_{n,2k+1,i}=x\bar{L}_{n,2k+1,i-1}-\frac{x}{1-x}\bar{L}_{n-1,2k-1,i+l-1}(x^2y,1/x)-
%\frac{x}{1-x}\bar{L}_{n-1,2k-1,i+l}(x^2y,1/x).
%\]
%\hfill $\Box$

We will combine the previous results to get
a recurrence for $\bar{L}_{n,2k}(x,y)$.

\begin{proposition}
\begin{equation}
\bar{L}_{n,2k}= \left\{ 
\begin{array}{ll}
\bar{L}_{n-1,2k}(x^2y,1/x)+
\frac{1}{1-x}\bar{L}_{n-1,2k-1}(x^2y,1/x) & \mbox{if $n > 2k$}\\
\frac{1}{1-x}\bar{L}_{n-1,2k-1}(x^2y,1/x) & \mbox{if $n = 2k$}
\end{array} \right.
\end{equation}
\label{L2krecurr}
\end{proposition}
\noindent{\bf Proof.} 
The $n=2k$ case is (\ref{Lnrec}).  Assume $n>2k$.
From Proposition \ref{decompose},
for $i=l(n-2k)+r$ with $2\le r\le n-2k$, we have
\[
\bar{L}_{n,2k+1,i}=x\bar{L}_{n,2k+1,i-1}-
\frac{x}{1-x}\bar{L}_{n-1,2k-1,i+l}(x^2y,1/x).
\]
and for $i=l(n-2k)+1$,
\[
\bar{L}_{n,2k+1,i}=x\bar{L}_{n,2k+1,i-1}-\frac{x}{1-x}\bar{L}_{n-1,2k-1,i+l-1}(x^2y,1/x)-
\frac{x}{1-x}\bar{L}_{n-1,2k-1,i+l}(x^2y,1/x).
\]
We sum
\begin{equation*}
 \bar{L}_{n,2k+1}  = 
\sum_{i=1}^{\infty} \bar{L}_{n,2k+1,i}  = 
 x\bar{L}_{n,2k}+x\sum_{i=1}^{\infty} \bar{L}_{n,2k+1,i}
-\frac{x}{1-x}\sum_{i=1}^{\infty} \bar{L}_{n-1,2k-1,i}(x^2y,1/x).
\end{equation*}
Therefore
\[
\bar{L}_{n,2k+1}=x\bar{L}_{n,2k}+x\bar{L}_{n,2k+1}
-\frac{x}{1-x}\bar{L}_{n-1,2k-1}(x^2y,1/x).
\]
Using Proposition \ref{map1},  that
\[
(1-x)\bar{L}_{n,2k+1}=x\bar{L}_{n-1,2k}(x^2y,1/x),
\]
we get the result. \hfill $\Box$

Now we can compute the generating function we seek.

\begin{theorem}  For $n \geq k \geq 0$,
\[
\bar{L}_{n,k}(x,y)=\frac{(x^{\lfloor k/2\rfloor+1}y^{\lfloor k/2\rfloor})
^{\lceil k/2 \rceil} \left[\begin{matrix}
n-\lceil k/2 \rceil\\ \lfloor k/2\rfloor
\end{matrix}\right]_{xy}}{(x;xy)_{\lceil k/2 \rceil}
(x^ny^{n-1};(xy)^{-1})_{\lfloor k/2\rfloor}}.
\]
\label{oetlh=}
\end{theorem}
{\bf Proof.} For $k=0$ the result holds,
%LLL fixed bar and took out k=1 case
%that is
%\[
%\bar{L}_{n,0}=1\ \ bar{L}_{n,1}=\frac{x}{1-x}.
%\]
as $\bar{L}_{n,0}=1$.
Let $n \geq m > 0$ and assume inductively that
the theorem is true for $(n,m-1)$, $(n-1,m-1)$,
and, if $n>m$, for $(n-1,m)$.

For the odd case $m=2k+1$, by Proposition \ref{map1},
\[
\bar{L}_{n,2k+1}=\frac{x}{1-x}\bar{L}_{n-1,2k}(x^2y,x^{-1}).
\]
By the induction hypothesis,
\[
\bar{L}_{n-1,2k}=\frac{(x^{k+1}y^k)^k\left[\begin{matrix}
n-1-k\\ k\end{matrix}\right]_{xy}}{(x;xy)_k(x^{n-1}y^{n-2};(xy)^{-1})_k},
\]
and substituting in the previous equation gives the result.\\

For the even case $m=2k$, we have by Proposition \ref{L2krecurr},
\begin{equation}
\bar{L}_{n,2k}= \left\{
\begin{array}{ll}
\bar{L}_{n-1,2k}(x^2y,1/x)+
\frac{1}{1-x}\bar{L}_{n-1,2k-1}(x^2y,1/x) & \mbox{if $n > 2k$}\\
\frac{1}{1-x}\bar{L}_{n-1,2k-1}(x^2y,1/x) & \mbox{if $n = 2k$}.\\
\end{array}
\right.
\end{equation}
%\[
%\bar{L}_{n,2k}=\bar{L}_{n-1,2k}(x^2y,1/x)+
%\frac{1}{1-x}\bar{L}_{n-1,2k-1}(x^2y,1/x).
%\]
By the induction hypothesis,
\[
\bar{L}_{n-1,2k-1}=\frac{(x^{k}y^{k-1})^k\left[\begin{matrix}
n-1-k\\ k-1\end{matrix}\right]_{xy}}{(x;xy)_k(x^{n-1}y^{n-2};(xy)^{-1})_{k-1}}.
\]
When $n=2k$, substituting this in the previous equation gives the result.
When $n>2k$, by the induction hypothesis, we also have
\[
\bar{L}_{n-1,2k}=\frac{(x^{k+1}y^k)^k\left[\begin{matrix}
n-1-k\\ k\end{matrix}\right]_{xy}}{(x;xy)_k(x^{n-1}y^{n-2};(xy)^{-1})_k}.
\]
%Therefore
%LLL  In denominator, last n should be n-1;  last xy should be inverse
%\[
%\bar{L}_{n-1,2k}(x^2y,1/x)=\frac{(x^{k+2}y^{k+1})^{k}\left[\begin{matrix}
%n-1-k\\ k\end{matrix}\right]_{xy}}{(x^2y;xy)_k(x^{n}y^{n};xy)_k},
%\]
%\[
%\bar{L}_{n-1,2k}(x^2y,1/x)=\frac{(x^{k+2}y^{k+1})^{k}\left[\begin{matrix}
%n-1-k\\ k\end{matrix}\right]_{xy}}{(x^2y;xy)_k(x^{n}y^{n-1};(xy)^{-1})_k},
%\ \ \ \
%\bar{L}_{n-1,2k-1}(x^2y,1/x)=\frac{(x^{k+1}y^{k})^k\left[\begin{matrix}
%n-1-k\\ k-1\end{matrix}\right]_{xy}}{(x^2y;xy)_k(x^{n}y^{n};(xy)^{-1})_{k-1}},
%\]
Then let
\[
\bar{L}_{n,2k} =
\bar{L}_{n-1,2k}(x^2y,1/x)+
\frac{1}{1-x}\bar{L}_{n-1,2k-1}(x^2y,1/x)=\frac{N_{n,k}}{D_{n,k}}\]
with
%LLL last n should be n-1
%\[
%D_{n,k}=(x;xy)_{k+1}(x^{n}y^{n};(xy)^{-1})_k.
%\]
\[
D_{n,k}=(x;xy)_{k+1}(x^{n}y^{n-1};(xy)^{-1})_k.
\]
Then
\[
N_{n,k}=(x^{k+2}y^{k+1})^{k}\left[\begin{matrix}
n-1-k\\ k\end{matrix}\right]_{xy}(1-x)+(x^{k+1}y^{k})^k\left[\begin{matrix}
n-1-k\\ k-1\end{matrix}\right]_{xy}(1-x^{n-k+1}y^{n-k}).
\]
We use these identities to simplify the numerator.
\begin{equation}
\left[\begin{matrix} n\\k \end{matrix}\right]_{q}=
\left[\begin{matrix} n-1\\k \end{matrix}\right]_{q}+
q^{n-k}\left[\begin{matrix} n-1\\k-1 \end{matrix}\right]_{q};
\ \ \ \ \ \ \ \  \ \ \ \ 
\left[\begin{matrix} n\\k \end{matrix}\right]_{q}=
\left[\begin{matrix} n-1\\k-1 \end{matrix}\right]_{q}+
q^{k}\left[\begin{matrix} n-1\\k \end{matrix}\right]_{q}  
\label{qbidents}
\end{equation}
So,
\begin{eqnarray*}
N_{n,k}&=&(x^{k+1}y^{k})^k \left( \left[\begin{matrix}
n-1-k\\ k-1\end{matrix}\right]_{xy}+x^ky^k \left[\begin{matrix}
n-1-k\\ k\end{matrix}\right]_{xy}\right)\\
& & -
(x^{k+1}y^{k})^k x^{k+1}y^{k} \left(\left[\begin{matrix}
n-1-k\\ k\end{matrix}\right]_{xy}+x^{n-2k}y^{n-2k}\left[\begin{matrix}
n-1-k\\ k-1\end{matrix}\right]_{xy}\right)\\
&=& (x^{k+1}y^{k})^k (1-x^{k+1}y^{k})\left[\begin{matrix}
n-k\\ k\end{matrix}\right]_{xy}.\ \ \ \ \ \
\end{eqnarray*}
\hfill $\Box$

The Odd/Even Truncated Lecture Hall Theorem (\ref{tlhxy}) is an
immediate consequence of Theorem \ref{oetlh=}.

\subsection{Truncated Anti-Lecture Hall Compositions}

Recall that for $n\ge k-1$, ${A}_{n,k}$ is the set
of compositions 
such that ~:
\[
\frac{\lambda_1}{n-k+1}\ge\ldots  \ge\frac{\lambda_k}{n}\ge 0.
\]
Our goal is to compute the generating function
\[
A_{n,k}(x,y)=\sum_{\lambda\in A_{n,k}}x^{|\lambda_o|}y^{|\lambda_e|}.
\]

We will again adapt the mapping of \cite{BME1}.
For $n \geq k \geq 1$, define
\[
\Theta_{n,k}: {A}_{n,k-1}\times \mathbb{N} \rightarrow {A}_{n,k},
\]
by $\Theta_{n,k}(\lambda,s) = (\mu_1, \mu_2, \ldots , \mu_k)=\mu$, where
\begin{eqnarray*}
\mu_1&\leftarrow &\left\lceil \frac{(n-k+1)\lambda_1}{n-k+2}\right\rceil+s \\
\mu_{2\ell} & \leftarrow &\lambda_{2\ell-1}, \ \ 1\le \ell\le k/2;\\
\mu_{2\ell+1}&\leftarrow &\left\lceil \frac{(n+2\ell-k+1)\lambda_{2\ell+1}}
{n+2\ell-k+2}\right\rceil+
\left\lfloor \frac{(n+2\ell-k+1)\lambda_{2\ell-1}}{n+2\ell-k}\right\rfloor
-\lambda_{2\ell},
\ \ 1\le \ell\le (k-1)/2. 
\end{eqnarray*}
Similar to the mapping
${\rm BME}_{n,k}$ of the previous section, it can be checked that
$\Theta_{n,k}$ is one-to-one and
$\Theta_{n,k}(A_{n,k-1} \times \mathbb{N}) \subseteq A_{n,k}$.
Furthermore, when $k$ is odd, $\Theta_{n,k}$ is onto $ A_{n,k}$ and
if $\mu = \Theta_{n,k}(\lambda,s)$, 
%This map will send a partition $\lambda$ in 
%${A}_{n,2k}$ and an integer $s$ 
%into a partition $\mu$ in ${A}_{n,2k+1}$~:
%\begin{eqnarray*}
%\mu_1&\leftarrow &\lceil \frac{(n-2k)\lambda_1}{n-2k+1}\rceil+s \\
%\mu_{2\ell} & \leftarrow &\lambda_{2\ell-1}; \ \ \ \ 1\le \ell\le k\\
%\mu_{2\ell+1}&\leftarrow &\lceil \frac{(n-2(k-\ell))\lambda_{2\ell+1}}{n-2(k-\ell)+1}\rceil+
%\lfloor \frac{(n-2(k-\ell)+2)\lambda_{2\ell-1}}{n-2(k-\ell)+1}\rfloor
%-\lambda_{2\ell};\ \ 1\le \ell\le k \\
%\end{eqnarray*}
%It is easy to check that $\mu\in A_{n,2k+1}$, 
then $|\mu_e|=|\lambda_o|$ and $|\mu_o|=2|\lambda_o|-|\lambda_e|+s$.
This implies~: 
\begin{proposition} For $n\ge 2k$,
\begin{equation}
A_{n,2k+1}(x,y) = 
\frac{1}{1-x}A_{n,2k}(x^2y,x^{-1}).
\end{equation}
\label{map2}
\end{proposition}
Again the recurrence for even $k$ is more difficult.
For $i\geq 0$,
let $A_{n,k,i}$ be the set of compositions in $A_{n,k}$ with $k^{th}$ part equal
to $i$.
It can be checked that $\Theta_{n,2k}$ maps 
${A}_{n,2k-1,i}\times \mathbb{N} $ bijectively to
${A}_{n,2k,i}$.  If $\mu = \Theta_{n,2k}(\lambda,s)$ it is not too hard
to see that
%\begin{eqnarray*}
%\mu_1&\leftarrow &\lceil \frac{(n-2k+1)\lambda_1}{n-2k+2}\rceil+s \\
%\mu_{2\ell} & \leftarrow &\lambda_{2\ell-1}; \ \ \ \ 1\le \ell\le k\\
%$\mu_{2\ell+1}&\leftarrow &\lceil \frac{(n-2(k-\ell)+1)\lambda_{2\ell+1}}{n-2(k-\ell)+2}\rceil+
%\lfloor \frac{(n-2(k-\ell)+3)\lambda_{2\ell-1}}{n-2(k-\ell)+2}\rfloor
%-\lambda_{2\ell};\ \ 1\le \ell\le k-1 \\
%\end{eqnarray*}
%
$\mu\in A_{n,2k+1}$, 
$|\mu_e|=|\lambda_o|$ and $|\mu_o|=2|\lambda_o|-|\lambda_e|+s-i-l$ with
$l=\lfloor i/n\rfloor$.
This implies~: 
\begin{proposition} For $n\ge 2k$,
\begin{equation}
A_{n,2k,i}(x,y) = 
\frac{x^{-i-l}}{1-x}A_{n,2k-1,i}(x^2y,x^{-1}).
\end{equation}
with $l=\lfloor i/n\rfloor$.
\label{bmeb1}
\end{proposition}
\begin{proposition} For $n\ge 2k$,
\begin{equation}
(1-x)A_{n,2k+1}=A_{n-1,2k}-\frac{x}{1-x}A_{n-1,2k-1}(x^2y,1/x),
\end{equation}
\label{oddeq}
\end{proposition}
\noindent{\bf Proof.} Let $i=ln+r$ with $0\le r\le n-1$.
It follows from the definitions that $A_{n,k,0}=A_{n-1,k-1}$. 
Consider the case $r=0$ and $l>0$.
If $\lambda \in A_{n,2k+1,i-1}$, then $\lambda_{2k}/(n-1) \geq (i-1)/n$,
so $\lambda_{2k}/(n-1) \geq i/n$.  Therefore adding 1 to the last part
of $\lambda$ gives a composition in $A_{n,2k+1,i}$.  Thus,
$ A_{n,2k+1,i}=xA_{n,2k+1,i-1}$.
For $r>0$, we have
\[
A_{n,2k+1,i}=xA_{n,2k+1,i-1}-x^iA_{n-1,2k,i-l-1}.
\]
To see this,
adding one to the $(2k+1)^{th}$ part of a composition
in $A_{n,2k+1,i-1}$ gives a composition in $A_{n,2k+1,i}$
unless the $(2k)^{th}$ part was equal to $i-l-1$.

By Proposition \ref{bmeb1}, 
 we get that
\[
A_{n-1,2k,i-l-1}=\frac{x^{1-i}}{1-x}A_{n-1,2k-1,i-l-1}(x^2y,1/x).
\]
We apply this and get
if $r\ge 1$,
\[
A_{n,2k+1,i}=xA_{n,2k+1,i-1}-\frac{x}{1-x}A_{n-1,2k-1,i-l-1}(x^2y,1/x);
\]
if $r=0$ and $l>0$
\[
A_{n,2k+1,i}=xA_{n,2k+1,i-1};
\]
otherwise, $i=0$ and 
$A_{n,k,0}=A_{n-1,k-1}$.
Now we sum on $i$~: 
\[
\sum_{i=1}^\infty A_{n,2k+1,i}=x\sum_{i=1}^\infty A_{n,2k+1,i-1}
-\frac{x}{1-x}\sum_{r=1}^{n-1}\sum_{l=0}^\infty
A_{n-1,2k-1,ln+r-l-1}(x^2y,1/x).
\]
This gives 
\[
A_{n,2k+1}-A_{n-1,2k}=xA_{n,2k+1}-\frac{x}{1-x}A_{n-1,2k-1}(x^2y,1/x),
\]
and therefore the result.
\hfill 
$\Box$

Now we get a recurrence for the number of even parts.
\begin{proposition}
\[
A_{n,2k}=\left\{ \begin{array}{ll} \frac{1}{1-x}A_{2k-1,2k-1}(y,x) & {\rm if}\ 
n=2k-1\\
A_{n-1,2k}(1/y,xy^2)
+\frac{y}{1-y}A_{n-1,2k-1} &{\rm if}\ n\ge 2k.\end{array}\right. 
\]
\label{A2krecurr}
\end{proposition}
\noindent{\bf Proof.} 
For the case $n=2k-1$. 
These are the objects such that
\[
\frac{\lambda_1}{0}\ge \frac{\lambda_2}{1}\ge \ldots \ge 
\frac{\lambda_{2k}}{2k-1}\ge 0.
\]
As the division by 0 makes the first inequality always valid if $\lambda_1\ge 0$, then
\[
A_{2k-1,k}(x,y)=\frac{1}{1-x}A_{2k-1,2k-1}(y,x).
\]

For $n\ge 2k$, we use the previous Proposition \ref{oddeq}.
From Proposition \ref{map2}, we know that 
$(1-x)A_{n,2k+1}=A_{n,2k}(x^2y,1/x)$. Therefore,
\[
A_{n,2k}(x^2y,1/x)=A_{n-1,2k}-\frac{x}{1-x}A_{n-1,2k-1}(x^2y,1/x).
\]
We make the substitutions $x=1/y$ and $y=xy^2$ and get the result.
\hfill $\Box$

\noindent{\bf Remark.} Note that, using Proposition \ref{A2krecurr}, we get~:
\[
A_{2k,2k}=
\frac{1}{1-y}(A_{2k-1,2k-1}-yA_{2k-1,2k-1}(xy^2,1/y)). 
\]

\begin{theorem}
{\rm (The Odd/Even Truncated Anti-Lecture Hall Theorem)}
For $n\ge k-1$ and $n,k\ge 0$,
\[
A_{n,k}(x,y)=\frac{\left[\begin{matrix}n\\ \lfloor k/2\rfloor \end{matrix}\right]_{xy}}
{(x;xy)_{\lceil k/2\rceil }(x^{n-k+1}y^{n-k+2};xy)_{\lfloor k/2\rfloor}}.
\]
\label{oetalh}
\end{theorem}

\noindent
{\bf Proof.}
We will prove the result by induction.
We know that
$A_{n,0}=1$.
Let $n+1 \geq m >0$  and assume inductively that
the theorem is true for $(n,m-1)$, $(n-1,m-1)$,
and, if $n\ge m$, for $(n-1,m)$.
For the odd case  $m=2k+1$,  by the induction hypothesis,
\[
A_{n,2k}=\frac{\left[\begin{matrix}n\\k\end{matrix}\right]_{xy}}
{(x;xy)_k(x^{n-2k+1}y^{n-2k+2};xy)_k}.
\]
Apply Proposition \ref{map2},
\[
A_{n,2k+1}=\frac{1}{1-x}A_{n,2k}(x^2y,1/x),
\]
and get the result.

For the even case $m=2k$, if $n=2k-1$, 
by the induction hypothesis, using Proposition \ref{A2krecurr}~:
\[
A_{2k-1,2k}(x,y)=\frac{1}{1-x}A_{2k-1,2k-1}(y,x)=
\frac{\left[\begin{matrix}2k-1\\k-1\end{matrix}\right]_{xy}}
{(y;xy)_k(x;xy)_k}.
\]

Now for $n\ge 2k$, by the induction hypothesis, 
using Proposition \ref{A2krecurr}~:
\[
A_{n-1,2k}(1/y,xy^2)=\frac{\left[\begin{matrix}n-1\\k\end{matrix}\right]_{xy}}
{(1/y;xy)_k(x^{n-2k+1}y^{n-2k+2};xy)_k}; \ \ \ \
A_{n-1,2k-1}=\frac{\left[\begin{matrix}n-1\\k-1\end{matrix}\right]_{xy}}
{(x;xy)_k(x^{n-2k+1}y^{n-2k+2};xy)_{k-1}}.
\]
\[
A_{n-1,2k}(1/y,xy^2)
-\frac{1/y}{1-1/y}A_{n-1,2k-1}
=\frac{N_{n,k}}{D_{n,k}}
\]
with $D_{n,k}=(1/y;xy)_{k+1}(x^{n-2k+1}y^{n-2k+2};xy)_k$.
Then
\begin{eqnarray*}
N_{n,k}&=&\left[\begin{matrix}n-1\\k\end{matrix}\right]_{xy}(1-x^ky^{k-1})
-1/y 
\left[\begin{matrix}n-1\\k-1\end{matrix}\right]_{xy}
(1-x^{n-k}y^{n-k+1})\\
&=&\left[\begin{matrix}n-1\\k\end{matrix}\right]_{xy}+
x^{n-k}y^{n-k}\left[\begin{matrix}n-1\\k-1\end{matrix}\right]_{xy}-\frac{1}{y}
\left(x^ky^k\left[\begin{matrix}n-1\\k\end{matrix}\right]_{xy}
+\left[\begin{matrix}n-1\\k-1\end{matrix}\right]_{xy}\right)\\
&=&(1-1/y)\left[\begin{matrix}n\\k\end{matrix}\right]_{xy},
\end{eqnarray*}
where the last step follows from (\ref{qbidents}).
\hfill $\Box$

\section{Combinatorial Characterizations and Refinements}

\subsection{Characterization of Truncated Lecture Hall Partitions}

In this section we characterize truncated lecture hall partitions
in terms of partitions into odd parts with certain restrictions.
 We first need a few steps.
\begin{proposition}
Given $n,t,j,l$, the weight generating function for 
the number of partitions into odd parts
in $\{1,3,...,2n-1\}$
in which exactly $t$ of the parts can be chosen from the set
$\{2j+1, 2j+3, ..., 2l-1\}$ is
\[
\frac{
q^{t(2j+1)}}
{(q;q^2)_{j}(q^{2n-1};q^{-2})_{n-l}}
\left[\begin{matrix}l-j-1+t\\ t\end{matrix}\right]_{q^2}.
\]
\label{restrictlem}
\end{proposition}
\noindent
{\bf Proof.}
The fraction, ignoring the denominator,
is the generating function for the 
partitions into the
odd parts which are not restricted.
Now count partitions into exactly $t$ parts from the
set $\{2j+1, 2j+3, ..., 2l-1\}$.
Take off $2j+1$ from each part. This is counted by $q^{t(2j+1)}$.
We are left with  a partition into even parts in a $t \times
2(l-j-1)$ box. \hfill $\Box$

\begin{corollary}
The generating function for the number of partitions into odd parts
in which exactly $t$ of the parts are greater than or equal to $2j+1$ is
\[
\frac{q^{t(2j+1)}}{(q;q^2)_j(q^2;q^2)_t}.
\]
\label{cor}
\end{corollary}

Applying Proposition \ref{restrictlem}  gives
the following.
\begin{proposition}
Let $R_{n,k}$
be  the set of partitions into odd parts less than  or equal to $2n-1$ 
where at most  $\lfloor k/2\rfloor$
parts can be chosen from the set
$
\{2\lceil k/2\rceil +1,2\lceil k/2\rceil +3,\ldots ,
2(n-\lfloor k/2\rfloor)-1\}$.
\begin{equation}
R_{n,k}(q) \triangleq
\sum_{\lambda \in R_{n,k}} q^{\lambda}=
\frac{
\sum_{i=0}^{\lfloor k/2 \rfloor}q^{i(2\lceil k/2 \rceil+1)}
\left[\begin{matrix}n-k-1+i\\ i\end{matrix}\right]_{q^2}}
{(q;q^2)_{\lceil k/2 \rceil}(q^{2n-1};q^{-2})_{\lfloor k/2 \rfloor}}
\label{toprove}.
\end{equation}
\hfill $\Box$
\label{Rnk}
\end{proposition}
Let $\bar{R}_{n,k}=R_{n,k}-R_{n,k-1}$ and let  $\bar{R}_{n,k}(q)$ be the 
corresponding generating function.
Recall that $\bar{L}_{n,k}$ is the set of lecture hall partitions in 
$L_{n,k}$ with $k$ positive parts.
%Let $L_{n,k}(q)=L_{n,k}(1,1,q)$ and 
%$\bar{L}_{n,k}(q)=\bar{L}_{n,k}(1,1,q)$. Then for $k>0$,
\begin{proposition}  For $k>0$,
\[
\bar{L}_{n,k}(q) \triangleq
 \sum_{\lambda \in \bar{L}{n,k}} q^{|\lambda|}
 = 
\frac{
q^{k+1\choose 2}\left[\begin{matrix}n-\lceil k/2 \rceil\\\lfloor k/2 \rfloor
\end{matrix}\right]_{q^2}}
{(q;q^2)_{\lceil k/2 \rceil}(q^{2n-1};q^{-2})_{\lfloor k/2 \rfloor}}
\]
\label{diff}
\end{proposition}
\noindent{\bf Proof.} Let $x=y=q$ in Theorem \ref{oetlh=}.
 %From Theorem \ref{rtlh=}, with $u=v=1$, we get
%\[
%\bar{L}_{n,k}(q)= 
%q^{k+1\choose 2}\left[\begin{matrix}n\\ k\end{matrix}\right]_q
%\frac{(-q^{n-k+1};q)_k}{(q^{2n-k+1};q)_k},
%\]
%and
%\begin{eqnarray*}
%\left[\begin{matrix}n\\ k\end{matrix}\right]_q
%\frac{(-q^{n-k+1};q)_k}{(q^{2n-k+1};q)_k}
%&  = &\frac{(q^{n-k+1};q)_k(-q^{n-k+1};q)_k}{(q;q)_k(q^{2n-k+1};q)_k}\\
%&  = &\frac{(q^{2(n-k+1)};q^2)_k} {(q;q^2)_{\lceil k/2\rceil}
%(q^2;q^{2})_{\lfloor k/2 \rfloor}(q^{2n-k+1};q)_k}\\
%&=& \frac{(q^{2(n-k+1)};q^2)_{\lfloor k/2 \rfloor}}
%{(q^2;q^{2})_{\lfloor k/2 \rfloor}}
%\frac{(q^{2n-2\lceil k/2 \rceil +2};q^2)_{\lceil k/2\rceil}}{(q;q^2)_{\lceil k/2\rceil}
%(q^{2n-k+1};q)_k}\\
%&=& \left[\begin{matrix}n-\lceil k/2 \rceil\\\lfloor k/2 \rfloor
%\end{matrix}\right]_{q^2}\frac{1}
%{(q;q^2)_{\lceil k/2 \rceil}(q^{2n-1},q^{-2})_{\lfloor k/2 \rfloor}}
%\end{eqnarray*}
\hfill $\Box$

\begin{theorem}
{\rm (Characterization of Truncated Lecture Hall Partitions)}
The number of truncated lecture partitions
of $N$ in $L_{n,k}$ is equal to the number of partitions of $N$ into
odd parts less than $2n$, with the following constraint on the parts:
at most
$\lfloor k/2\rfloor$
parts can be chosen from the set
$$
\{2\lceil k/2\rceil +1,2\lceil k/2\rceil +3,\ldots ,
2(n-\lfloor k/2\rfloor)-1\}.$$
\label{combi}
\end{theorem}
\noindent{\bf Proof.} 
We must show that
$L_{n,k}(q)\triangleq \sum_{m=0}^k \bar{L}_{n,m} = R_{n,k}(q)$.
Since  $L_{n,0}(q)= R_{n,0}(q)=1$,
it suffices to show that for  $k>0$,
$\bar{L}_{n,k}(q)=
\bar{R}_{n,k}(q).$
We use $\bar{L}_{n,k}(q)$ from Proposition \ref{diff}, and now compute 
$\bar{R}_{n,k}(q)={R}_{n,k}(q)-{R}_{n,k-1}(q).$

For $k>0$, first look at the case $k=2s$.  Using Proposition \ref{Rnk},
\begin{equation*}
\bar{R}_{n,2s}(q) =
\frac{
\sum_{i=0}^{s}q^{i(2s+1)}
\left[\begin{matrix}n-2s-1+i\\ i\end{matrix}\right]_{q^2} -
(1-q^{2(n-s)+1})\sum_{i=0}^{s-1}q^{i(2s+1)}
\left[\begin{matrix}n-2s+i\\ i\end{matrix}\right]_{q^2}
}
{(q;q^2)_{s}(q^{2n-1};q^{-2})_{s}}.
\end{equation*}
Since the denominators of $\bar{R}_{n,2s}(q)$ and
$\bar{L}_{n,2s}(q)$ agree, we focus on the numerator $E(n,s)$~:
 $$E(n,s)=\sum_{i=0}^{s}q^{i(2s+1)}
\left[\begin{matrix}n-2s-1+i\\ i\end{matrix}\right]_{q^2} -
(1-q^{2(n-s)+1})\sum_{i=0}^{s-1}q^{i(2s+1)}
\left[\begin{matrix}n-2s+i\\ i\end{matrix}\right]_{q^2}.$$
%We will use the fact that
%~:
%\[
%\left[\begin{matrix} n\\k \end{matrix}\right]_{q}=
%\left[\begin{matrix} n-1\\k \end{matrix}\right]_{q}+
%q^{n-k}\left[\begin{matrix} n-1\\k-1 \end{matrix}\right]_{q}.
%\]

Therefore, using the first identity of (\ref{qbidents}),
\begin{eqnarray*}
\sum_{i=0}^{s-1}q^{i(2s+1)}
\left(\left[\begin{matrix}n-2s-1+i\\ i\end{matrix}\right]_{q^2}
-\left[\begin{matrix}n-2s+i\\ i\end{matrix}\right]_{q^2}\right) 
& = &
-q^{2(n-2s)}\sum_{i=0}^{s-1}q^{i(2s+1)}
\left[\begin{matrix}n-2s+i-1\\ i-1\end{matrix}\right]_{q^2}\\
& = &
-q^{2(n-s)+1}\sum_{i=0}^{s-2}q^{i(2s+1)}
\left[\begin{matrix}n-2s+i\\ i\end{matrix}\right]_{q^2}.
\end{eqnarray*}
We then get~: 
\begin{eqnarray*}
E(n,s)&=&q^{s(2s+1)}
\left[\begin{matrix}n-s-1\\ s\end{matrix}\right]_{q^2}
+q^{2(n-s)+1} q^{(s-1)(2s+1)}
\left[\begin{matrix}n-s-1\\ s-1\end{matrix}\right]_{q^2}\\
&=&q^{s(2s+1)}\left[\begin{matrix}n-s\\ s\end{matrix}\right]_{q^2}
= q^{k+1\choose 2}\left[\begin{matrix}n-\lceil k/2 \rceil\\ \lfloor k/2 \rfloor\end{matrix}\right]_{q^2}.
\end{eqnarray*}
Now when $k=2s+1$, applying Proposition \ref{Rnk} gives
\begin{equation*}
\bar{R}_{n,2s+1}(q) =
\frac{
\sum_{i=0}^{s}q^{i(2s+3)}
\left[\begin{matrix}n-2s-2+i\\ i\end{matrix}\right]_{q^2} -
(1-q^{2s+1})\sum_{i=0}^{s}q^{i(2s+1)}
\left[\begin{matrix}n-2s-1+i\\ i\end{matrix}\right]_{q^2}
}
{(q;q^2)_{s+1}(q^{2n-1};q^{-2})_{s}}.
\end{equation*}
Again we can focus on the numerator $O(n,s)$~:
\[
O(n,s)=
\sum_{i=0}^{s}q^{i(2s+3)}
\left[\begin{matrix}n-2s-2+i\\ i\end{matrix}\right]_{q^2} -
(1-q^{2s+1})\sum_{i=0}^{s}q^{i(2s+1)}
\left[\begin{matrix}n-2s-1+i\\ i\end{matrix}\right]_{q^2}.
\]
Then, using the second identity of (\ref{qbidents}),
%\[
%\left[\begin{matrix} n+k\\k \end{matrix}\right]_{q^2}=
%\left[\begin{matrix} n+k-1\\k-1 \end{matrix}\right]_{q^2}+
%q^{2k}\left[\begin{matrix} n+k-1\\k \end{matrix}\right]_{q^2}.
%\]
%Then
\begin{eqnarray*}
\sum_{i=0}^{s}q^{i(2s+1)}
\left(q^{2i}\left[\begin{matrix}n-2s-2+i\\ i\end{matrix}\right]_{q^2} -
\left[\begin{matrix}n-2s-1+i\\ i\end{matrix}\right]_{q^2}\right)&=&
-\sum_{i=1}^{s}q^{i(2s+1)}
\left[\begin{matrix}n-2s-2+i\\ i-1\end{matrix}\right]_{q^2}\\
&=&
-q^{2s+1}\sum_{i=0}^{s-1}q^{i(2s+1)}
\left[\begin{matrix}n-2s-1+i\\ i\end{matrix}\right]_{q^2}.
\end{eqnarray*}
Applying this to $O(n,s)$ gives the result ~:
\begin{eqnarray*}
O(n,s)&=&q^{2s+1}q^{s(2s+1)}\left[\begin{matrix}n-s-1\\ s\end{matrix}\right]_{q^2}
= q^{k+1\choose 2}\left[\begin{matrix}n-\lceil k/2\rceil\\ \lfloor k/2\rfloor\end{matrix}\right]_{q^2}.
\end{eqnarray*}
$\Box$

\subsection{Finitizations of Refinements of Euler's Theorem}

Euler's Theorem says that 
the number of partitions of $N$ into distinct parts
is equal to the number of partitions of $N$ into odd parts.
(A straightforward proof shows $(-q;q)_\infty=1/(q;q^2)_\infty$.)
The Lecture Hall Theorem is a finitization of that theorem~:
the number of partitions of $N$ in $L_n$
is equal to the number of partitions of $N$ into odd parts
less than $2n$. 
We will show how the truncated lecture hall results
give finitizations of refinements of Euler's theorem
that can be easily from (slight modifications of) Sylvester's bijection.

\noindent
{\bf Finitization 1.}
\begin{equation}
\sum_{m=0}^kq^{m+1\choose 2}\left[\begin{matrix}n\\ m\end{matrix}\right]_q
\frac{(q^{n-m+1};q)_m}{(q^{2n-m+1};q)_m}
=
\sum_{i=0}^{\lfloor k/2 \rfloor}
\frac{
q^{i(2\lceil k/2 \rceil+1)}
\left[\begin{matrix}n-k-1+i\\ i\end{matrix}\right]_{q^2}}
{(q;q^2)_{\lceil k/2 \rceil}(q^{2n-1};q^{-2})_{\lfloor k/2 \rfloor}}.
\label{ident1}
\end{equation}

This is 
Theorem \ref{combi}, which states that $L_{n,k}(q)=R_{n,k}(q)$, where 
$L_{n,k}(q)$ and $R_{n,k}(q)$ are given by
Theorem \ref{rtlh} and Proposition \ref{Rnk}, respectively.

Taking limits as $n \rightarrow \infty$ in (\ref{ident1}) gives
\[
\sum_{m=0}^k \frac{q^{m+1\choose 2}}
{(q;q)_m}
=
\sum_{i=0}^{\lfloor k/2 \rfloor}
\frac{
q^{i(2\lceil k/2 \rceil+1)}}
{(q;q^2)_{\lceil k/2 \rceil}(q^2;q^2)_{i}}.
\]
Note that $q^{m+1\choose 2}/(q;q)_m$ is the generating function for
partitions into $m$ distinct parts.  Applying Corollary \ref{cor} to the
right hand side, then, gives the following refinement of Euler's Theorem.
\begin{quote}
{\bf Refinement 1:}
The number of partitions of $N$ into at most $k$ distinct parts
is equal to the number of partitions of $N$ into odd parts such
that at most $\lfloor k/2 \rfloor$ of the parts are greater than
or equal to
$2 \lceil k/2 \rceil +1$.
\end{quote}
\begin{quote}
Theorem \ref{combi} itself is a further refinement.
\end{quote}

\noindent {\bf Finitization 2.}
\begin{equation}
q^{k+1\choose 2}\left[\begin{matrix}n\\ k\end{matrix}\right]_q
\frac{(q^{n-k+1};q)_k}{(q^{2n-k+1};q)_k}=
\frac{
q^{\binom{k+1}{2}}
\left[\begin{matrix}n-\lceil k/2 \rceil\\\lfloor k/2 \rfloor
\end{matrix}\right]_{q^2}}
{(q;q^2)_{\lceil k/2 \rceil}(q^{2n-1};q^{-2})_{\lfloor k/2 \rfloor}}.
\label{ident2}
\end{equation}
This is Proposition \ref{diff} combined with eq. (\ref{rtlh=}),
setting $u=v=1$ in eq. (\ref{rtlh=}).

Letting $n \rightarrow \infty$ in (\ref{ident2}) gives
\begin{equation}
\frac{q^{\binom{k+1}{2}}}
{(q;q^2)_k}=
\frac{q^{\binom{k+1}{2}}}
{(q;q^2)_{\lceil k/2 \rceil}(q^{2};q^{-2})_{\lfloor k/2 \rfloor}}.
\label{ident2lim}
\end{equation}
On the right-hand side of (\ref{ident2lim}),
if we use
\[q^{\binom{k+1}{2}}=
q^{\lceil k/2 \rceil-\lfloor k/2 \rfloor}
q^{\lfloor k/2 \rfloor(2 \lceil k/2 \rceil +1)}
\]
and apply Corollary \ref{cor}, (\ref{ident2lim})
can be read as the following refinement of 
Euler's theorem.

\begin{quote}
{\bf Refinement 2:}
The number of partitions of $N$ into exactly $k$ distinct parts
is equal to the number of partitions of
$N-( \lceil k/2 \rceil-\lfloor k/2 \rfloor)$ into odd parts such
that  exactly $\lfloor k/2 \rfloor$ of the parts are greater than
or equal to
$2 \lceil k/2 \rceil +1$.
\end{quote}

A combinatorial interpretation of (\ref{ident2})
 will give a further refinement
of Refinement 2.  To get this, note that the right-hand side of
(\ref{ident2}) can be written as
\[\left(\frac{q^{\lceil k/2 \rceil-\lfloor k/2 \rfloor}}
{1-q^{2n-2\lfloor k/2 \rfloor +1}}\right)
\frac{
q^{\lfloor k/2 \rfloor(2 \lceil k/2 \rceil +1)}
\left[\begin{matrix}n-\lceil k/2 \rceil\\\lfloor k/2 \rfloor
\end{matrix}\right]_{q^2}}
{(q;q^2)_{\lceil k/2 \rceil}(q^{2n-1};q^{-2})_{\lfloor k/2 \rfloor-1}}.
\]
Using Proposition \ref{restrictlem}, the last quotient above is the generating
function for the number of partitions into odd parts less than
$2n$ with exactly $\lfloor k/2 \rfloor$ parts in the set
$\{2\lfloor k/2 \rfloor+1, \ldots, 2n-2\lfloor k/2 \rfloor+1\}$.
So (\ref{ident2})  says:
\begin{quote}
{\bf Further refinement 2:}
 The number of partitions of $N$ in $\bar{L}_{n,k}$
is equal to the number of partitions of
$N-( \lceil k/2 \rceil-\lfloor k/2 \rfloor)$ into odd parts less than
$2n$ with at least $\lfloor k/2 \rfloor$  parts
in $\{2\lfloor k/2 \rfloor+1, \ldots, 2n-2\lfloor k/2 \rfloor+1\}$,
but at most
$\lfloor k/2\rfloor$ parts
less than $2n-2\lfloor k/2\rfloor+1$.
\end{quote}

\noindent {\bf Finitization 3.}

Using Theorem \ref{oetlh=} with $x=zq$ and $y=q/z$, we get
\[
\bar{L}_{n,k}(zq,q/z)=\frac{zq^{2\lfloor k/2\rfloor+1})
^{\lceil k/2 \rceil} \left[\begin{matrix}
n-\lceil k/2 \rceil\\ \lfloor k/2\rfloor
\end{matrix}\right]_{q^2}}{(zq;q^2)_{\lceil k/2 \rceil}
(zq^{2n-1};(q^2)^{-1})_{\lfloor k/2\rfloor}}.
\]
Now let $n\rightarrow \infty$ and ${\mathcal D}_k$ be the set
of partitions into $k$  distinct parts  and get
\[
\sum_{\lambda\in {\mathcal D}_{2k-1}\cup {\mathcal D}_{2k}}
z^{|\lambda_o|-{|\lambda_e|}}q^{|\lambda|}
=\frac{z^k q^{k(2k-1)}}{(zq;q^2)(q^2;q^2)_k}.
\]

We say that the Durfee rectangle size of a partition $\lambda$ into odd parts
is $k$ if  $\lambda_k\ \ge 2k-1$ and $\lambda_{k+1}\le 2k-1$.
 
\begin{quote}
{\bf Refinement 3:}
The number of partitions $\lambda$ 
of $N$ into  $2k-1$ or $2k$ distinct parts
with $|\lambda_o|-|\lambda_e|=j$
is equal to the number of partitions of $N$ into $j$ parts and
Durfee rectangle size $k$.
\end{quote}

\subsection{Interpretation of Anti-Lecture Hall Theorems}

We consider the limiting case for truncated anti-lecture hall
compositions.  Note that for fixed $k$, as $n \rightarrow \infty$,
the set $A_{n,k}$ approaches $P_k$, the set of partitions into 
$k$ nonnegative parts.  From Theorem \ref{oetalh} we have
\[
\sum_{\lambda \in A_{n,k}} x^{|\lambda_o|}y^{|\lambda_e|} =
\frac{\left[\begin{matrix}n\\ \lfloor k/2\rfloor \end{matrix}\right
]_{xy}}
{(x;xy)_{\lceil k/2\rceil }(x^{n-k+1}y^{n-k+2};xy)_{\lfloor k/2\rfloor}}.
\]
Taking limits as $n \rightarrow \infty$ and substituing $x=zq$ and $y=q/z$, 
gives
\[
\sum_{\lambda \in P_{k}} q^{|\lambda|}z^{|\lambda_o|-|\lambda_e|}
= \frac{1}{(zq;q^2)_{\lceil k/2\rceil} (q^2;q^2)_{\lfloor k/2\rfloor}}, 
\]
which is the well known ``transpose theorem'':  The number of partitions 
$\lambda$ of $N$ into $k$ nonnegative parts and 
$|\lambda_o|-|\lambda_e|=j$ is equal to 
the number of partitions of
$N$ with largest part less than or equal to $k$ and $j$ of the parts odd.

Similarly, consider the identity which combines
Theorem \ref{talhuv}, setting $u=v=1$ with
Theorem \ref{oetalh}, setting $x=y=q$:
\begin{equation}
\left[\begin{matrix}n\\ k\end{matrix}\right]_q
\frac{(-q^{n-k+1};q)_k}{(q^{2(n-k+1)};q)_k} =
\frac{\left[\begin{matrix}n\\ \lfloor k/2\rfloor \end{matrix}\right
]_{q^2}}
{(q;q^2)_{\lceil k/2\rceil }(q^{2n-2k+3};q^2)_{\lfloor k/2\rfloor}}.
\label{Aident}
\end{equation}
Taking limits as $n \rightarrow \infty$ shows that (\ref{Aident})
 can be interpreted
as a finitization of the transpose theorem.

\section{Conclusion}

We hopefully have demonstrated that basic hypergeometric $q$-series
are a good tool to give refined Lecture Hall-Type
Theorem proofs. 
Furthermore, the BME mapping plays a significant role
in the development of 2-variable generating functions for anti-Lecture
Hall compositions and truncated objects.
Our study of the 2-variable generating function
of the truncated objects leads to the following $x,y$-series identity.
Using Equations (\ref{lhpxy}) and (\ref{tlhxy}), we have
\begin{proposition}
\[
\sum_{m=0}^n\left(x^{\lfloor m/2\rfloor+1}y^{\lfloor m/2\rfloor}\right)^{\lceil m/2\rceil}\frac{\left[\begin{matrix} n-\lceil m/2\rceil\\ \lfloor m/2\rfloor\end{matrix}\right]_{xy}}{(x;xy)_{\lceil m/2\rceil}(x^ny^{n-1};(xy)^{-1})_{\lfloor m/2\rfloor}}=\prod_{i=1}^n \frac{1}
{1-x^iy^{i-1}}.
\]
\end{proposition}

As a special case, substituing $x=a$ and $y=q/a$, we get the identity~:
\begin{proposition}
\[
\sum_{m=0}^n\left(aq^{\lfloor m/2\rfloor}\right)^{\lceil m/2\rceil}\frac{\left[\begin{matrix} n-\lceil m/2\rceil\\ \lfloor m/2\rfloor\end{matrix}\right]_{q}}{(a;q)_{\lceil m/2\rceil}(aq^{n-1};q^{-1})_{\lfloor m/2\rfloor}}=1/(a;q)_n.
\]
\end{proposition}
In future work we will consider 
truncated versions of the $(k,l)$-lecture hall partitions
of \cite{BME2}.

\bibliography{bib13}
\bibliographystyle{plain}

\noindent{\bf Acknowledgments.} The authors want to thank the 
Universit\'e de Versailles Saint-Quentin and in particular the laboratory
PRiSM for inviting Dr Savage for the month of December 2002. They also
want to thank the CNRS and NSF for their joint-grant.
\end{document}